\documentclass[article]{siamart220329}
\usepackage{amsmath,amssymb}
\usepackage{comment}
\usepackage{url}

\graphicspath{{./figures/}}

\usepackage{algorithm}
\usepackage{algpseudocode}
\usepackage{array,booktabs}

\newdimen{\algindent}
\algnewcommand{\LeftComment}[1]{\Statex \hspace{\algindent}
\( \triangleright \) #1}

\usepackage{subfigure} 
\usepackage{threeparttable} 

\def\x{{\bf x}}
\def\y{{\bf y}}
\def\dt{{\Delta t}}
\def\I{\mathcal{I}}
\def\V{\mathcal{V}}
\def\n{{\bf n}}
\def\ub{{\bf u}}

\newcommand{\be}{\begin{equation}}
\newcommand{\ee}{\end{equation}}
\newcommand{\ba}{\begin{aligned}}
\newcommand{\ea}{\end{aligned}}

\def\ccm{Center for Computational Mathematics, Flatiron Institute, Simons Foundation,
  New York, New York 10010}

\def\umich{Department of Mathematics, University of Michigan,
Ann Arbor, Michigan 48109\\}

\def\nyu{Courant Institute, New York University, New York, New York 10012}

\def\thu{Yau Mathematical Sciences Center,
  Tsinghua University, Beijing China 100084\\}

\def\papertitle{Space-time adaptive methods for parabolic evolution equations}

\title{\papertitle}
\author{
  Jun Wang\thanks{\thu\,
    ({\tt jwang2020@tsinghua.edu.cn}).} \and
  Jie Su\thanks{\thu\,
    ({\tt suj22@mails.tsinghua.edu.cn}).} \and
Leslie Greengard\thanks{\ccm\ and\ \nyu \,
    ({\tt lgreengard@flatironinstitute.org}).} \and
Shidong Jiang\thanks{\ccm\,
    ({\tt sjiang@flatironinstitute.org}).} \and
Shravan Veerapaneni\thanks{\umich
     ({\tt shravan@umich.edu}).}
}

\begin{document}

\maketitle

\begin{abstract}
We present a family of integral equation-based solvers 
for the heat equation, reaction-diffusion systems, the unsteady Stokes
equation and the
incompressible Navier-Stokes equations in two space dimensions.
Our emphasis is on the development of methods
that can efficiently follow complex
solution features in \emph{space-time} by refinement and coarsening at each time step
on an adaptive quadtree.
For simplicity, we focus on problems posed in a square 
domain with periodic boundary conditions.
The performance and robustness of the methods are illustrated with several 
numerical examples. 
\end{abstract}



\pagestyle{myheadings}
\thispagestyle{plain}
\markboth{J. Wang, J. Su, L. Greengard, S. Jiang, S. Veerapaneni}
{Space-time adaptive methods for diffusion-dominated problems}

\section{Introduction \label{sec:intro}}  
Many problems in scientific and
engineering applications 
require the solution of scalar diffusion equations
of the form
\begin{equation} \label{heatfree}
\begin{split}
u_{t}(\x,t) &= D \Delta u(\x,t)+F(u,\nabla u,\x,t), \\
u(\x,0) &= u_0(\x),
\end{split}
\end{equation}
for $t>0$ and $\x \in B \subset \mathbb{R}^d$ 
or systems of the form 
\begin{equation} \label{heatsys}
\begin{split}
\ub_{t}(\x,t) &= {\bf D} \Delta \ub(\x,t)+{\bf F}(\ub,\nabla \ub,\x,t), \\
\ub(\x,0) &= \ub_0(\x),
\end{split}
\end{equation}
where ${\bf D}\in \mathbb{R}^{p\times p}$ is a diagonal matrix of diffusion constants and
$\ub(\x,t) \in \mathbb{R}^p$.
Here, $u_0(\x)$ and $\ub_0(\x)$ are used to denote the initial data, $F(\cdot)$ and ${\bf F}(\cdot)$ are
given forcing functions.
For simplicity, we focus on the two dimensional 
case  ($d=2$) and assume the solution satisfies periodic boundary 
conditions on $\Omega = [-\frac{1}{2},\frac{1}{2}]^2$:
\begin{equation}
\begin{aligned}
\ub(-\tfrac{1}{2},y) &= \ub(\tfrac{1}{2},y), \quad
\ub_x(-\tfrac{1}{2},y) = 
\ub_x(\tfrac{1}{2},y),\ {\rm for}\ 
-\tfrac{1}{2} \leq y \leq \tfrac{1}{2}, \\
\ub(x,-\tfrac{1}{2}) &= \ub(x,\tfrac{1}{2}), \quad
\ub_y(x,-\tfrac{1}{2}) = 
\ub_y(x,\tfrac{1}{2}),\ {\rm for}\ 
-\tfrac{1}{2} \leq x \leq \tfrac{1}{2}.
\end{aligned}
\end{equation}
As a matter of nomenclature, we refer to the  
general equation \eqref{heatsys} 
as {\em semilinear} if $\mathbf{F}$ depends on $\ub$ but not
its gradient, and as {\em linear} when $\mathbf{F} = \mathbf{F}(\x,t)$.
A longer introduction to parabolic potential theory for problems of the type
considered here was posted on arXiv \cite{wang2019arxiv} in 2019, including a discussion
of boundary value problems in moving domains. Some of the material below is drawn
from that earlier presentation.

Classical methods, such as finite difference or finite element methods,
evolve the solution at successive time steps, using some local 
approximation of the partial differential equation (PDE)
\eqref{heatsys} itself.
While quite general in applicability,
they are subject to a severe constraint on the time step when
relying on explicit marching schemes. Using the forward Euler method,
for example, if we approximate the Laplacian using the second order
central difference operator on a uniform grid, then 
the stability condition for solving the scalar initial value problem requires
that 
$\Delta t \leq \frac{h^2}{2D}$, where $\Delta t$ is the time 
step, and $h = \Delta x = \Delta y$ is the step in the spatial discretization.
In general, the stability constraint takes the form
$\dt=\mathcal{O}(h^2/D)$ where $h$ is the {\em smallest} spacing in the
spatial discretization.
This constraint is typically overcome
by using an implicit marching scheme. 
In the linear case, this
requires the solution of a large linear system of equations
at each time step, coupling all the spatial grid points. 
In the semilinear case, a fully implicit discretization scheme
requires the solution of a large, nonlinear system of equations.
The available literature on such approaches is vast and we 
do not seek to review it here.
Instead, we focus on the use of parabolic potential theory
\cite{costabel,friedman1964,guenther1988,pogorzelski} and 
recently developed fast algorithms that will permit us to solve 
such initial value problems on fully adaptive meshes with relative
ease.

{\sc Remark 1.1}
{\em
Global spectral methods based on the fast 
Fourier transform (FFT) are, of course, powerful
alternatives for solving the problems considered here on a periodic
box. However, they evolve the solution on a uniform grid 
\cite{boyd2001,shen_spectral,trefethen_spectral} and are unsuitable for 
space-time adaptivity.
}

\subsection{Linear problems}\label{sec:linear}

Let us first consider the linear version of the scalar problem
\eqref{heatfree}, where the forcing term $F(u,\nabla u,\x,t) = F(\x,t)$
is known. This equation is well-posed
\cite{friedman1964,guenther1988,pogorzelski}, 
with the solution given 
at any later time $t$ in closed form as
\begin{equation}
u(\x,t) = \I[u_0](\x,t)+\V[F](\x,t),
\label{freesol}
\end{equation}
with 
\begin{equation}
\I[u_0](\x,t) = \int_{\Omega} G(\x-\y,t) u_0(\y) \, d\y
\label{initpot}
\end{equation}
and 
\begin{equation}
\V[F](\x,t) = \int_0^t\int_{\Omega} 
G(\x-\y,t-\tau) F(\y,\tau) \, d\y d\tau,
\label{volpot}
\end{equation}
where
\begin{equation}
G(\x,t) =  \sum_{\n \in \mathbf{Z}^2} 
\frac{e^{-\|\x - \n \|^2/4Dt}}{(4 \pi D t)} \ 
 =  \sum_{\n \in \mathbf{Z}^2} 
e^{-4 \pi^2 |\n|^2 D t} 
e^{2 \pi i \n \cdot \x} 
\label{heatker}
\end{equation}
is the periodic Green's function for the heat equation 
(the {\em heat kernel}) in two dimensions,
expressed using either the method of images or as a Fourier series.
The functions $\I[u_0]$ and $\V[F]$ are referred to as {\em initial} and
{\em volume} heat potentials, respectively.

A compelling feature of the representation 
\eqref{freesol} is that it is exact, explicit, and does not involve the inversion
of any matrices. It requires only the 
{\em evaluation} of the initial
and volume potentials. Thus, stability never arises as an issue, and the error
is simply the quadrature error in approximating the integral operators
that define the relevant potentials themselves.
This approach is typically ignored for one simple reason:
in the absence
of fast algorithms, the computational cost is excessive. More
precisely, assuming there are $N_T$ time steps,
and $N$ points in the discretization of the domain, naive summation
applied to \eqref{volpot} requires $O(N^2 N_T^2)$ work.
It is also the case that the heat kernel is singular as a function of time so that
accurate and efficient quadrature rules are required, but those issues are
well-understood.

It is, however, straightforward to see that, for a time step $\dt$, 
\begin{equation}
\begin{aligned}
u(\x,t+\dt) = 
&\int_{\Omega} G(\x-\y,\dt) \, 
u(\y,t) \, d\y + \\
&\qquad \qquad \qquad \int_{t}^{t+\dt}
\int_{\Omega} G(\x-\y,t+\dt-\tau) \, F(\y,\tau) \, d\y d\tau.
\end{aligned}
\label{vmarch}
\end{equation}
That is to say, the heat equation can be solved by 
convolving the heat kernel with the current solution at every time step,
followed by the addition of a {\em local} contribution from 
the forcing term $F(\x,t)$ over the last time step alone.
This can be viewed as a marching scheme for the PDE expressed in 
integral form, and reflects the semigroup property of heat flow.
Using this approach, the apparent obstacle of
history-dependence vanishes and the net cost drops from
$O(N^2 N_T^2)$ to $O(N^2 N_T)$ work, assuming some suitable discretization/quadrature
rule has been applied to the integrals in \eqref{vmarch}. Fast algorithms,
discussed below, will reduce this further
to $O(N N_T)$ work.
Linear systems are treated in the same manner: each component of \eqref{heatsys}
advances as in \eqref{vmarch} with its own diffusion constant.

\subsection{Parabolic systems of equations}
For general systems of the form  \eqref{heatsys}, 
the solution can still be expressed using the vector version of \eqref{freesol}, 
with the solution given 
at time $t$ as
\begin{equation}
\ub(\x,t) = \I[\ub_0](\x,t)+\V[\bf F,\ub](\x,t)
\label{freesolsys}
\end{equation}
where
\begin{equation}
\I[\ub_0](\x,t) = \int_{\Omega} G(\x-\y,t) \ub_0(\y) \, d\y, 
\label{initpotsys}
\end{equation}
and 
\begin{equation}
\V[{\bf F}, \ub](\x,t) = \int_0^t\int_{\Omega} 
G(\x-\y,t-\tau) {\bf F}(\ub,\nabla \ub,\y,\tau) \, d\y d\tau.
\label{volpotsys}
\end{equation}

There are significant advantages to this approach when compared with elliptic marching methods,
discussed in section \ref{sec:AMmethod}.

\subsection{The unsteady Stokes equations}\label{sec:us}
Linearization of the incompressible Navier-Stokes equations leads to the
unsteady Stokes equations:
\be\label{stokes}
\begin{split}
    \mathbf{u}_t(\mathbf{x},t) &= \nu\nabla^2 \mathbf{u} - \bigtriangledown p +\mathbf{F}(\mathbf{x},t) ,\quad \nabla\cdot\mathbf{u} = 0 \\   
    \mathbf{u}(\mathbf{x},0) &= \mathbf{u}_0(\mathbf{x}).
    \end{split}
\ee
Here, $\mathbf{u}(\mathbf{x},t)$ is the velocity field of interest, $\nu$ is the viscosity, and $p(\mathbf{x}, t) $ is the pressure.
We assume the forcing term {\bf F}(\x,t) is given.
To apply the scheme used for the systems of diffusion equations \eqref{heatsys}, we introduce the Helmholtz decomposition of the forcing term:

\be
\mathbf{F}(\mathbf{x}, t) = \mathbf{F}_S(\mathbf{x}, t) + \mathbf{F}_G(\mathbf{x}, t),
\ee
where \(\mathbf{F}_S\) and \(\mathbf{F}_G\) are the solenoidal and gradient components of $\mathbf{F}$, respectively.

It is well-known that, in free space, there is a simple construction for the Helmholtz decomposition of $\mathbf{F}$
\cite{ladyzhenskaya1964mathematical}. Namely, if $\mathbf{F} \in L^2(\mathbb{R}^d)$, then

\begin{eqnarray}
    \mathbf{F}_G &=&  \nabla \left( \nabla \cdot \int_{\mathbb{R}^d} G_L(\mathbf{x} - \mathbf{y}) \mathbf{F}(\mathbf{y}) \, d\mathbf{y} \right), \label{hdecomp1} \\
    \mathbf{F}_S &=& \mathbf{F} - \mathbf{F}_G,\label{hdecomp2} 
\end{eqnarray}
where $G_L(\mathbf{x}) = -\frac{1}{2\pi} \ln(|\mathbf{x}|)$ is the Green's function for the Laplace equation. Moreover,
the solution to the unsteady Stokes equations is given by
\cite{greengard2019new,ladyzhenskaya1964mathematical} 
\be\label{solunsteay}
\begin{split}
\mathbf{u}(\mathbf{x}, t) &= \mathcal{I}[\mathbf{u}_0](\mathbf{x}, t) + \mathcal{V}[\mathbf{F}_S](\mathbf{x}, t)\\
\nabla p(\mathbf{x}, t) &=  \mathbf{F}_G(\mathbf{x}, t)
\end{split}
\ee
The same holds for smooth, periodic vector fields in the box $D$, and
either the fast multipole method (FMM)
\cite{AskhamCerfon,cheng2006jcp,ethridge2001sisc,langston2011camcs,lee1996jcp,biros2015cicp,malhotra2016toms} or the dual-space multilevel kernel-splitting (DMK) method~\cite{jiang2025dual} can be used to efficiently compute 
\(\mathbf{F}_G\) and \(\mathbf{F}_S\) (see section \ref{sec:fmm}).

\subsection{The Navier-Stokes equations}\label{sec:ns}

In a periodic box, the incompressible Navier-Stokes equations take the form
\be\label{nstokes}
\begin{split}
    \mathbf{u}_t(\mathbf{x},t) &= \nu\nabla^2 \mathbf{u} - \bigtriangledown p -(\mathbf{u}\cdot \nabla)\mathbf{u} +\mathbf{f}(\mathbf{x},t) ,\quad \nabla\cdot\mathbf{u} = 0 \\   
    \mathbf{u}(\mathbf{x},0) &= \mathbf{u}_0(\mathbf{x}).
    \end{split}
\ee
From the preceding discussion, if we define 
\begin{equation} 
{\bf F}(\ub,\x,t) = - (\mathbf{u} \cdot \nabla) \mathbf{u} + \mathbf{f}(\mathbf{x}, t),
\end{equation}
(with $\mathbf{F} = 
\mathbf{F}_S+ \mathbf{F}_G$)
we may again represent the solution as \eqref{solunsteay}, where
\begin{equation}
\V[{\bf F}_S,\ub](\x,t) = 
\int_0^t\int_{\Omega} 
G(\x-\y,t-\tau) {\bf F}_S (\ub,\y,t) \, d\y d\tau.
\label{volpotsysn}
\end{equation}
Thus, the original PDE has been recast as the
nonlinear Volterra equation
\begin{equation}
\ub(\x,t) = \I[\ub_0](\x,t)+\V[{\bf F}_S, \ub](\x,t).
\end{equation}
Numerical methods for this are discussed in section \ref{sec:AMmethod}.

\subsection{Synopsis of the paper}

In section \ref{sec:prelim}, 
we briefly review the relevant fast algorithms for the evaluation of volume heat potentials,
the solution of the Poisson equation
and the Helmholtz decomposition of vector fields
\cite{ethridge2001sisc,fgt2024,lee1996jcp}.
Numerical schemes for linear and semilinear problems
are summarized in section \ref{sec:alg}, including automatic coarsening/refinement
strategies.
We illustrate the performance of our schemes in 
section \ref{sec:examples} and
discuss their extension to
domains with complicated boundary in section \ref{sec:conclusion}.

We hope to make clear some of the distinct advantages to be gained from using the integral equation framework.  These include the following:
\begin{enumerate}
\item 
Explicit marching schemes are 
{\em unconditionally stable} in the linear setting even
on highly nonuniform grids.
\item Fully implicit methods for scalar semilinear problems 
require only the solution of {\em scalar} nonlinear equations
(as noted previously in \cite{epperson2}). (PDE-based methods require the 
solution of nonlinear elliptic systems.)
\item
Linear scaling fast algorithms are available with arbitrary 
order accuracy in space and time.
\item
Automatic coarsening and refinement is straightforward and robust.
\end{enumerate}

\section{Mathematical Preliminaries\label{sec:prelim}}
Common to all of our adaptive algorithms is the need to automatically 
resolve a function $u(\x)$, which we 
assume to be smooth and periodic on the 
unit square $D$ but may have multiscale features.

\subsection{Adaptive resolution of a function} \label{sec:adap}

In order to resolve the
function
$u(\x)$ using an adaptive quadtree,
we assume that superimposed on $D$ is a hierarchy of refinements, with 
grid level 0 defined to be $D$ itself. Grid level $l+1$ is obtained from grid level $l$
by subdividing some squares at level $l$ into four equal parts.
For a square $B$ at level $l$, the four squares at level $l+1$ obtained by its 
subdivision will be referred to as its children. The refinement process continues
until a polynomial approximation of the form 
\begin{equation}
 u(\x) = u_0(x_1,x_2) \approx \sum_{n=0}^{K-1}\sum_{m=0}^{K-1}
u_{n,m} p_n(x_1) p_m(x_2)
\label{polyapprox}
\end{equation}
is determined to be sufficiently accurate on each leaf node of the tree. 
Here, $p_n(x)$ is the Chebyshev
polynomial $T_n(x)$ scaled to an interval of length $2^{-l}$ at level $l$. 

The determination of whether $u(\x)$ is resolved can be accomplished 
in several ways. Let $B$ be a leaf node with $u(\x)$ given on that box by
a Chebyshev expansion of the form \eqref{polyapprox}.
One can then evaluate $u(\x)$ on a $2K \times 2K$ grid
covering $B$ and compute the discrete
$L^2$ error, denoted by $E_2$, over these target points.
If $E_2 > \epsilon$, for the user-specified tolerance,
the leaf node $B$ is subdivided and the process repeated on each
of its children. Alternatively, one can assess the resolution 
through the error estimate \cite{trefethentail}
\begin{equation}
 E \approx \sqrt{ \sum_{n^2+m^2 \geq K^2} |u_{n,m}|^2}/K , 
\label{resolve_error}
\end{equation}
that is, by measuring the $L^2$ norm of the coefficients $\{ u_{n,m} \}$ outside 
the $L^2$ ball of radius $K$.
We will say that the tree obtained by a systematic use of such 
procedure starting with a single $K \times K$ grid on the unit box $D$
{\em resolves} $u(\x)$. 

We will assume, for the sake of simplicity, 
that the quadtree satisfies a fairly standard restriction, namely, that 
two leaf nodes which share a boundary point must be no more than one refinement level 
apart (Fig. \ref{fig:adaptree}). Trees satisfying this condition are generally called
{\em balanced} or {\em level-restricted}.
Once an adaptive tree is constructed to resolve $u(\x)$, 
however, it may not satisfy this criterion. However, assuming that 
the resolving tree so-constructed has $N_b$ leaf nodes and that its depth 
is of the order $O(\log N_b)$,
it is straightforward to balance the tree, enforcing the level-restriction,
in a subsequent sweep using $O(N_b \log N_b)$ time and storage.
We omit the details and refer the reader to 
\cite{quadtree,ethridge2001sisc,sundar2016}.
The number of degrees of freedom in such a discretization of $u(\x)$ is
$N = N_b \, K^2$.

\subsection{The continuous FGT} \label{sec:dgt}

\begin{figure}[htbp]
\centering
\includegraphics[width=\textwidth]{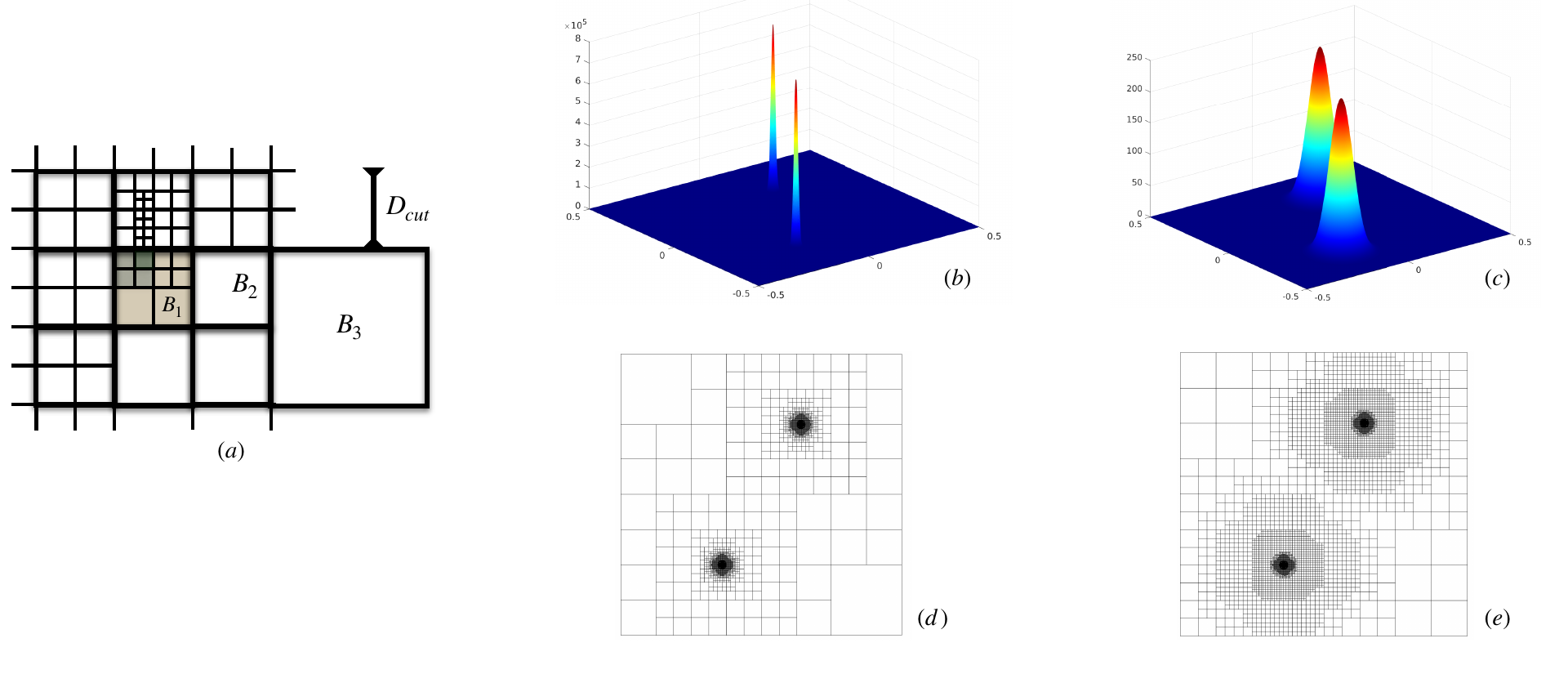}
\caption{
[Adapted from \cite{fgt2024}]
The new FGT computes
volume integrals discretized on a level-restricted
quadtree with arbitrary order accuracy, determined by the 
order of the polynomial approximation on each leaf node.
(a) For any precision $\epsilon$, the FGT defines 
the cutoff level $l$ in the tree hierarchy to be that where the box dimension 
$D_{cut}= 2^{-l}$ satisfies $e^{-D_{cut}^2/(4\delta)} < \epsilon$.
In the figure, $B_2$ is a leaf node at the cutoff level, $B_1$ is a leaf node at
a finer level and $B_3$ is a leaf node at a coarser level.
[For readers familiar with hierarchical fast transforms:
plane wave expansions are merged in an upward pass (from fine to coarse) until the 
cutoff level is reached. They are then translated to near neighbors at the cutoff
level and communicated to finer levels in a downward pass. For coarse leaf nodes 
such as $B_3$,
the continuous Gauss transform is evaluated in the near field {\em directly}, using
precomputed tables for optimal performance.]
An interesting aspect of the Gauss transform {\em on non-uniform data structures}
is that the grid resolving the output may need to be finer than the grid resolving 
the input! In (b) we consider $u_0(\x)$ as the sum of two sharply peaked Gaussians.
In (c), we plot the initial potential (the Gauss transform of $u_0$) for
$\delta = 4 \times 10^{-3}$.
In (d), we show 
the quadtree that resolves the {\em input} data to $12$ digits of accuracy.
In (e), we show
the quadtree that resolves the {\em output} data to $12$ digits of accuracy
after convolving with a sharp Gaussian.
The reason for this surprising behavior is that a sharp feature can diffused a 
short distance and reach regions where it is too sharp to be well resolved by the 
original adaptive grid, indicating that careful monitoring of resolution is
essential. 
The total number of boxes in (e) is actually slightly less than in (d),
after refinement where more resolution is needed and coarsening where
the resolution is determined to be excessive. See section \ref{coarseningsec} for further details.
\label{fig:adaptree}}
\end{figure}

The initial potential  \eqref{initpot} is an example of the 
{\em continuous Gauss transform}. 
Fast algorithms for the evaluation of such potentials
are known as fast Gauss transforms (FGTs) and have been available for 
several decades
\cite{brattkus1992siap,greengard2000acha,greengard1990cpam,greengard1991fgt,greengard1998nfgt,Lee2006,li2009sisc,sampath2010pfgt,spivak2010sisc,strain1994sisc,tausch2009sisc,veerapaneni2008jcp,wang2018sisc}.
In \cite{fgt2024}, a new version of the FGT was introduced that is particularly
efficient for continuous sources approximated as piecewise polynomials on adaptive tensor product grids with either free-space or periodic boundary conditions. This version of the FGT
exploits the separable structure of the Gaussian
to compute the potential (for any parameter $\delta$) at all grid points in $O(N)$  time with a small prefactor,
comparable in speed to the fast Fourier transform in work per grid point, even in its  fully adaptive form.
For a detailed description of the algorithm, see the original paper,
which is briefly summarized in Fig. \ref{fig:adaptree}.

\subsection{The Helmholtz decomposition of a vector field}
\label{sec:fmm} 

Suppose now that each component of the vector field $\mathbf{F}$ has been resolved on 
an adaptive level restricted tree in the manner described in 
section \ref{sec:adap}. Then, the Helmholtz decomposition \eqref{hdecomp1},
\eqref{hdecomp2} can 
be computed in $O(N) = O(N_b K^2)$ work using either DMK~\cite{jiang2025dual} or the FMM~\cite{ethridge2001sisc}. 
We refer the reader to the original papers for details. For our present purposes,
we simply note that these solvers work on exactly the same data structure as the FGT
above. 

{\sc Remark 2.1}
{\em
A distinct advantage of the volume integral formula in \eqref{hdecomp1} 
(with the integral over the box $D$ using the periodic Green's function) is that
the differential operators 
can be applied to the Green's function itself, rather than to the data.
Thus, the Helmholtz decomposition can be computed without loss of precision 
from numerical differentiation, as noted in the original 
FMM-based volume integral paper \cite{lee1996jcp}.
}

In this paper, we make use of the FMM-based Poisson solver 
\cite{ethridge2001sisc} which  computes the potential and its derivatives at 
all grid points in $O(N)$  time with a small prefactor,
comparable in speed to the fast Fourier transform in work per grid point, even in its 
fully adaptive form.

\section{Marching schemes for parabolic evolution equations} \label{sec:alg}

The marching form for parabolic equations is described above for the linear,
scalar case in \eqref{vmarch}. We restate it here in the case of systems with
nonlinear forcing terms. For this,
let us denote the $n$th time step by 
$t_n=n\dt$ and the solution at the $n$th time step by $\ub_n(\x)=\ub(\x,n\dt)$.
The exact evolution is given by
\be
\ub_{n+1}(\x)=\I[\ub_{n}](\x,\dt)+\int_{n\dt}^{(n+1)\dt}\mathbf{U}(\ub,\x,\tau)d\tau
\label{am1}
\ee
with
\be
\mathbf{U}(\ub,\x,\tau)=\int_{\Omega} G(\x-\y,t_{n+1}-\tau) \mathbf{F}(\ub(\y,\tau),\nabla \ub(\y,\tau),\y,\tau) \, d\y,
\label{am2}
\ee
where $G(\x,t)$ is given by \eqref{heatker}.
Recall that the problem is linear if 
\be 
\mathbf{F}(\ub(\y,\tau),\nabla \ub(\y,\tau),\y,\tau) =  
 \mathbf{F}(\y,\tau),
\ee
and semilinear (a reaction-diffusion system) if 
\be 
\mathbf{F}(\ub(\y,\tau),\nabla \ub(\y,\tau),\y,\tau) =  
 \mathbf{F}(\ub(\y,\tau),\y,\tau).
 \ee
With a slight abuse of notation, 
the formula \eqref{am2} also covers the unsteady Stokes equations and the Navier-Stokes equations, where
the vector field ${\bf F}$ should be understood as its solenoidal component (see \eqref{hdecomp2}): 
\be
\mathbf{U}(\ub,\x,\tau)=\int_{\Omega} G(\x-\y,t_{n+1}-\tau) \mathbf{F}_S(\ub(\y,\tau),\nabla \ub(\y,\tau),\y,\tau) \, d\y.
\label{am2a}
\ee
In the unsteady Stokes setting \eqref{solunsteay}, the evolution is linear and thus there is no dependence on $\ub$.

While a large variety of discretization schemes can be applied to \eqref{am1}, we will focus on
linear multistep methods. By analogy with marching schemes for ordinary 
differential equations (ODEs), we will refer to implicit schemes as Adams-Moulton methods and to 
explicit schemes as Adams-Bashforth methods.

\subsection{Adams-Moulton and Adams-Bashforth methods}\label{sec:AMmethod}

Linear multistep methods for \eqref{am1} are based on approximating each component of
the vector $\mathbf{U}(\x,\tau)$ for each $\x$ by an interpolating
polynomial in time, $p(\x,\tau)$, of degree $s-1$. For implicit marching, 
we include the current time in the interpolation and the 
conditions to be satisfied are
\be
p(\x,t_{n+1-i})=U(u_{n+1-i},\x,t_{n+1-i}), \qquad i=0,\ldots,s-1.
\label{am3}
\ee
For explicit marching, we use only previously known values and the
conditions to be satisfied are
\be
p(\x,t_{n+1-i})=U(u_{n+1-i},\x,t_{n+1-i}), \qquad i=1,\ldots,s.
\label{ab3}
\ee
Replacing the components of $\mathbf{U}$ in \eqref{am1} by their interpolants $p$ 
and integrating the resulting expression, we obtain
\be
\ub_{n+1}(\x)=\I[\ub_{n}](\x,\dt)+\dt \sum_{i=0}^{s} 
b_i \mathbf{U}_{n+1}(u_{n+1-i},\x,t_{n+1-i}),
\label{am4}
\ee
where $b_0=0$ for the explicit (Adams-Bashforth) schemes and
$b_s=0$ for the implicit (Adams-Moulton) schemes. Both achieve 
accuracy in time of order $s$, as in the case of ODEs, and the coefficients
$\{ b_i \}$ themselves are the same as for the classical integrators with the 
same name (see, for example, \cite{hwn_ode}), and listed 
for orders $1,2$ and $4$ in Table \ref{abmcoeffs}.

\renewcommand{\arraystretch}{1.2} 
\begin{table}[htbp]
 \caption{Parameters for linear multistep methods of various orders}
  \centering
  \begin{tabular}{l|l|l}
    \toprule
$s$ & Adams-Moulton & Adams-Bashforth \\
\midrule
1 & $b_0=1$  & $b_1=1$  \\
2 & $b_0 = \frac12,\ b_1 = \frac12$ & $b_1 = \frac32,\ b_2= -\frac12$    \\
4 & $b_0 = \frac9{24}, 
b_1 = \frac{19}{24}, b_2 = -\frac5{24}, b_3 = \frac1{24}$ &
$b_1 = \frac{55}{24}, b_2 = -\frac{59}{24}, b_3 = \frac{37}{24}, b_4 = -\frac{9}{24}$ \\ 
\bottomrule
\end{tabular}
\label{abmcoeffs}
\end{table}

For any such linear multistep method, at time $t_{n+1}$, 
the previously obtained solutions $u_{n+1-i}$ ($i=1,\dots,s$) can
be stored, so that the initial potential 
$\I[\ub_{n}](\x,\dt)$ and the spatial volume integrals
$\mathbf{U}_{n+1}(\ub_{n+1-i},\x,t_{n+1-i})$ for $i>0$ can be computed explicitly 
and in optimal time with the FGT.
Moreover, for Adams-Moulton methods, where $b_0 \neq 0$,
\be
\ba
\mathbf{U}_{n+1}(\ub_{n+1},\x,t_{n+1})&=
\lim_{\tau\rightarrow t_{n+1}}\int_{\Omega} G(\x-\y,t_{n+1}-\tau) \mathbf{F}(\ub(\y,\tau),\y,\tau) \, d\y\\
&=\mathbf{F}(\ub_{n+1},\x,t_{n+1}),
\label{am6}
\ea
\ee
using the $\delta$-function property of the Green's function.
For linear evolution equations, this completes the discussion of implicit marching. Note that
no equations need to be solved, unlike elliptic marching schemes.
For semilinear problems,
marching scheme of the form \eqref{am4}
leads to a set of {\it spatially uncoupled}
nonlinear equations:
\be
\ub_{n+1}(\x)-\dt b_0 \mathbf{F}(\ub_{n+1}(\x),\x,t_{n+1})=\mathbf{g}(\x),
\label{am7}
\ee
where
\be
\label{gdef}
\mathbf{g}(\x) = \I[\ub_{n}](\x,\dt)+ \dt \sum_{i=1}^{s-1} 
b_i \mathbf{U}_{n+1}(u_{n+1-i},\x,t_{n+1-i}).
\ee
This, again, implies that
solving \eqref{am4} is much easier than PDE-based marching which results in a nonlinear elliptic PDE to be
solved at each time step.  Moreover, there is no obstacle to obtaining high order accuracy in time.
In the PDE-based approach, it {\em is} possible to solve only scalar nonlinear equations 
through the use of {\em operator splitting} methods, but the order of accuracy
is then limited by the splitting error. While remedies for this are available
\cite{hansen,jahnke,mclachlan,murua,strang,tyson}, they are rather involved and
splitting methods are typically implemented with second order accuracy.

Many root finding methods can be used to solve \eqref{am7}. 
For the scalar case, it is convenient to use the (derivative-free)
secant method with two initial guesses chosen to be $u^{(0)}_{n+1}(\x)=u_{n}(\x)$
  and $u^{(1)}_{n+1}(\x)=g(\x)+\dt\, b_0\, F(u_{n},\x,t_{n})$. 
For the case of systems, with $\ub(\x,t) \in \mathbb{R}^p$, we use Newton's method at each
point $\x$ with initial guess
$\ub^{(0)}_{n+1}(\x)=\ub_{n}(\x)$.
For more general nonlinear forcing terms, where $\mathbf{F} = \mathbf{F}(\ub(\x,t),\nabla \ub(\x,t),\x,t)$,
such as the Navier-Stokes equations, the gradient term causes the nonlinear systems at time $t_{n+1}$
to be coupled and, in this paper, we will consider only explicit or preedictor-corrector methods.
(Thus, our approach to the Navier-Stokes equations is similar to that of standard PDE-based schemes 
which are typically implicit in the diffusion term and explicit in the forcing term.)

{\sc Remark 3.1}
{\em
  As with any high order linear multistep method, when the order of accuracy
  is greater than $2$, one needs an alternative scheme to compute the
  first few steps $\{u_1, \ldots, u_{s-1}\}$ with high order accuracy.
  In our current implementation, we do this by combining a second order accurate scheme
  with Richardson extrapolation to the desired order.
}

{\sc Remark 3.2}
{\em
For the Navier-Stokes equations, we use a predictor-corrector approach for 
each order of accuracy.
That is, we first compute an approximate solution 
$\ub^p_{n+1}$ at each new time step $t_{n+1}$ using an Adams-Bashforth method of order $s$. We then use 
$\ub^p_{n+1}$ in the Adams-Moulton formula of the same order to obtain an improved estimate for
$\ub_{n+1}$. This doesn't change the order of accuracy but the predictor-corrector scheme has a stability
region which includes part of the imaginary axis 
\cite{hwn_ode} and is more suitable for problems which involve advection. Significant improvements in stability can be obtained using the SAV scheme developed in \cite{shen2018sinum,shen2018jcp,shen2019sirev}. Recently, robust higher-order time-marching schemes have been constructed in \cite{huang2024sinum}.
}

\subsection{Spatial refinement and coarsening} \label{coarseningsec}

For the linear case, the quadtree structure at each time $t_n$ can be refined
in a straightforward manner, as discussed in section \ref{sec:adap}.
For the semilinear case, if refinement is needed, the value of the function $F(u,\x,t)$ 
must be obtained at a new target point $\x$ at time $t_{n+1}$. For this,
we solve \eqref{am7} to obtain $u(\x,t_{n+1})$ and
$F(u,\x,t_{n+1})$.
Note that the function $g(\x)$ in \eqref{am7} can be assumed to be resolved as it involves
the Gauss transform of previously computed solutions. 
An important feature of the Gauss transform is that it guarantees resolution of the output
on a Chebyshev grid at each leaf node of the adaptive data structure, so that interpolation 
is accurate at any new evaluation point. (This issue does not
arise for the Navier-Stokes equations, where we use explicit marching schemes, so that
refinement is activated only if 
either $u(\x,t_{n+1})$ or $F(u,\x,t_{n+1})$
is determined to be unresolved.

For coarsening, the reverse strategy is employed. Suppose $B$ 
is a box with four children.
We can compute $u(\x,t_{n+1})$ at the $K \times K$ Chebyshev nodes in $B$.
If the corresponding interpolant from $B$ matches both $u(\x,t_{n+1})$ and 
$F(u,\x,t_{n+1})$ to the desired precision at all child grid points, then we 
delete the child nodes and retain the coarsened grid on $B$ for subsequent time steps.
This is summarized in Algorithm~\ref{spatialadap}.

\begin{algorithm} 
\caption{Spatial refinement and coarsening algorithm for the semilinear heat equation}
\label{spatialadap}
\begin{algorithmic}
  \Procedure{spatial\_adap}{$T_{n}$, $F$, $l$, $\epsilon$, $T_{n+1}$, $u_n$}
  \LeftComment{Input: $T_{n}$ - the adaptive quadtree at $t_{n}$.\\
    $\qquad\qquad\qquad\, F$ - the function $F(u,\x,t)$.\\
    $\qquad\qquad\qquad\,\,\, l$ - the number of levels of $T_{n}$.\\
    $\qquad\qquad\qquad\,\,\, \epsilon$ - the desired precision.}
  \LeftComment{Output: $T_{n+1}$ - the adaptive quadtree at $t_{n+1}$.\\
    $\qquad\qquad\qquad\,\,\,\, u_{n+1}$ - the solution $u_{n+1}$ 
on $T_{n+1}$.}
        
  \LeftComment{{\bf Step 1} Refinement:}
  \For {each box $B$ in $T_{n}$}
  \If {$B$ is childless}
  \State Calculate the Chebyshev interpolant of $F(u_{n+1},\x,t_{n+1})$ on $B$ and
  compute $F(u_{n+1},\x,t_{n+1})$ on a refined grid corresponding to $B$'s children.
  Compute the error $E$ in the interpolant on the child grids.
  \If {$E < \epsilon$}
  \State \textbf{break}
  \Else
  \State Add the four child boxes to the quadtree data structure.
  \EndIf
  \EndIf
  \EndFor
  
  \LeftComment{{\bf Step 2} Coarsening:}
  \For {$i = l-1, 0, -1$}
  \State Calculate the values of the function $F(u_{n+1},\x,t_{n+1})$ 
on the $K\times K$ grid for each box that has 
  four child leaf nodes.
  \State Delete the four child boxes if $u_{n+1}$ and $F(u_{n+1},\x,t_{n+1})$ are resolved
  at the parent level.
  \EndFor
  
  \EndProcedure
\end{algorithmic}
\end{algorithm}

\section{Numerical examples}\label{sec:examples}

In this section, we illustrate the performance
of the methods described above for solving a variety of problems governed by
the parabolic evolution equations on the unit box $D= [-0.5,0.5]^2$. 
All of the component algorithms were implemented in Fortran
and all timings listed are generated using a laptop with a 2.3 GHz Intel Core i5 processor with 4GB RAM.
 
\subsection{The forced linear heat equation}

We first consider the behavior of our methods for the linear version of the
heat equation \eqref{heatfree} 
with periodic boundary conditions on the unit box, as discussed in section \ref{sec:linear}:
\begin{eqnarray*}
    u_t(\mathbf{x}, t) &=& \Delta u(\mathbf{x}, t) + F( \mathbf{x}, t), \nonumber \\
    u(\mathbf{x}, 0) &=& u_0(\mathbf{x}). 
\end{eqnarray*}

The complexity of the problem is controlled by the forcing function 
\( F(\mathbf{x}, t) \). 
To test our adaptive refinement/coarsening strategy, we let
\begin{equation}
    F(\mathbf{x},t) = \sum_{\mathbf{j}\in \mathbf{Z}^2} e^{-|\mathbf{x}-\mathbf{c}_1(t)-\mathbf{j}|^2/\delta} - 0.5e^{-|\mathbf{x}-\mathbf{c}_2(t)-\mathbf{j}|^2/\delta} \label{heat_forcing} 
\end{equation}
with $\mathbf{c}_1(t) = 0.25(\cos(20\pi t),\sin(20\pi t))$ and $\mathbf{c}_2(t) = 0.25(\cos(40\pi t+\pi),\sin(40\pi t+\pi))$.
The points  $\mathbf{c}_1,\mathbf{c}_2$ 
move on a circle at different velocities and the Gaussians are sharply peaked,
causing $F(\mathbf{x},t)$ to be highly inhomogeneous in space-time. 

\begin{figure}[htbp]
\centering
\subfigure[t=0.001, \#leaf boxes=592]{
\includegraphics[width=0.3\textwidth]{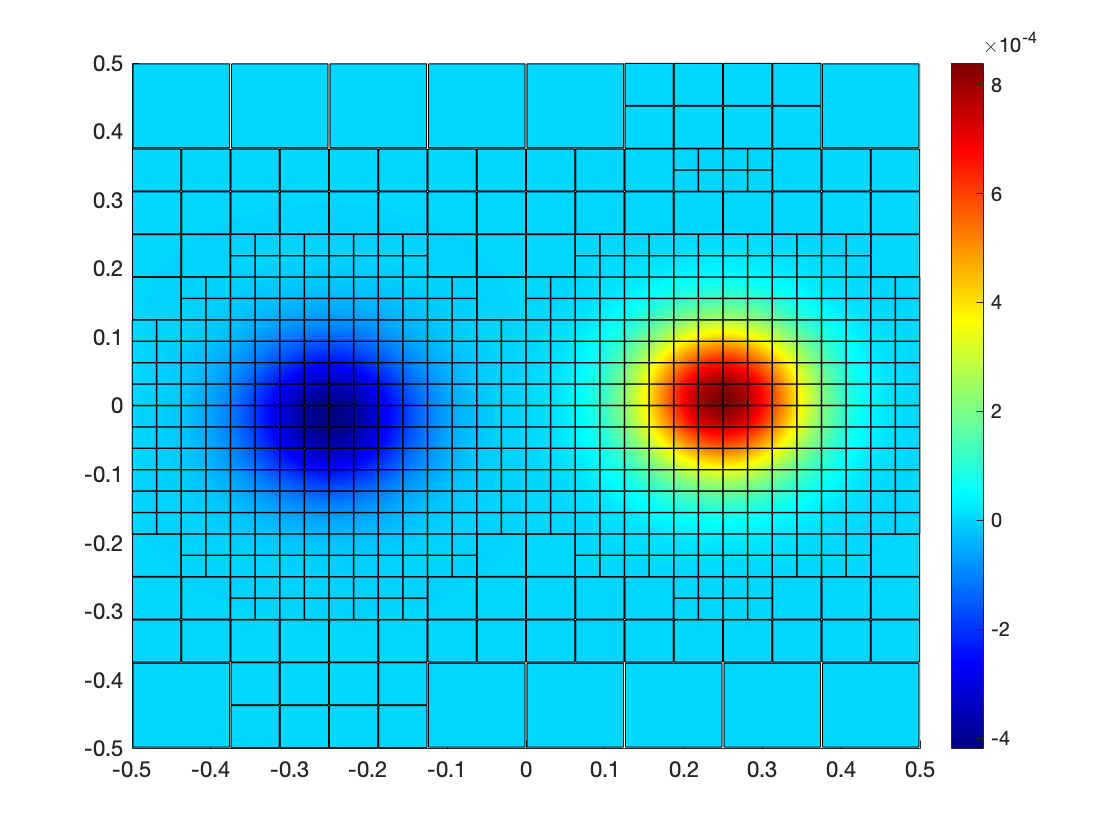}
}
\subfigure[t=0.02, \#leaf boxes=637]{
\includegraphics[width=0.3\textwidth]{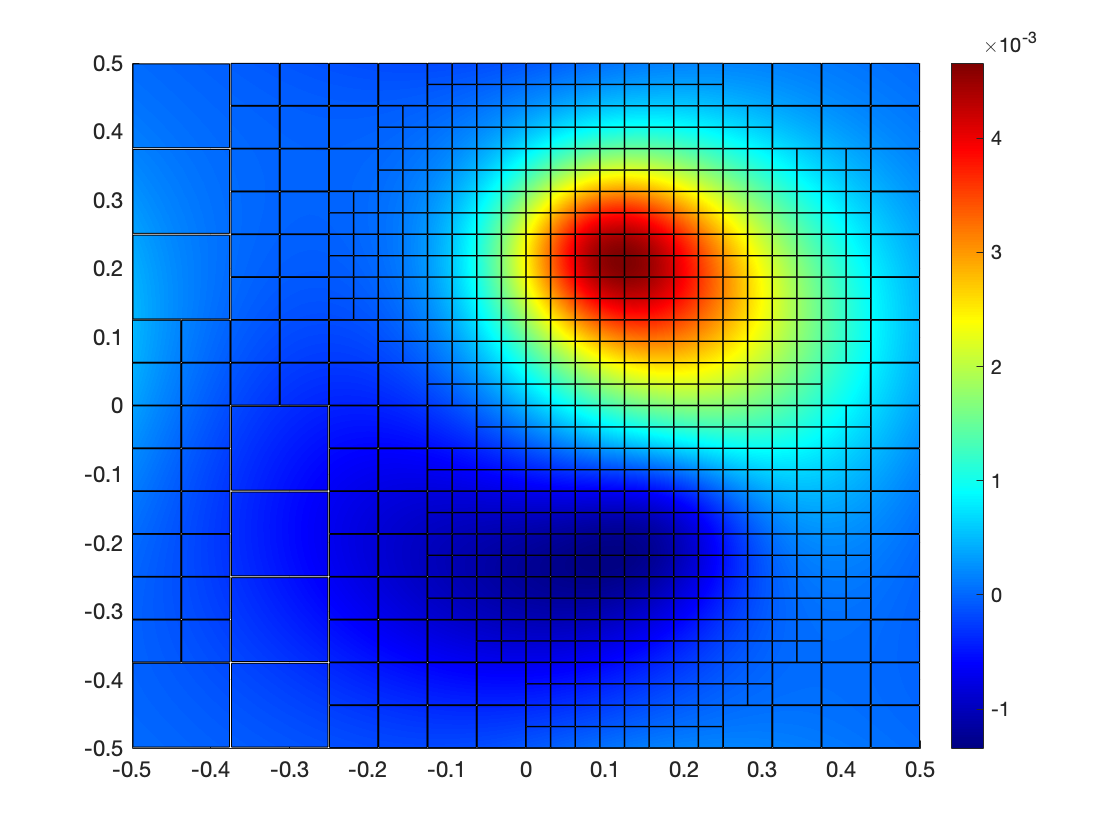}
}
\subfigure[t=0.05, \#leaf boxes=409]{
\includegraphics[width=0.3\textwidth]{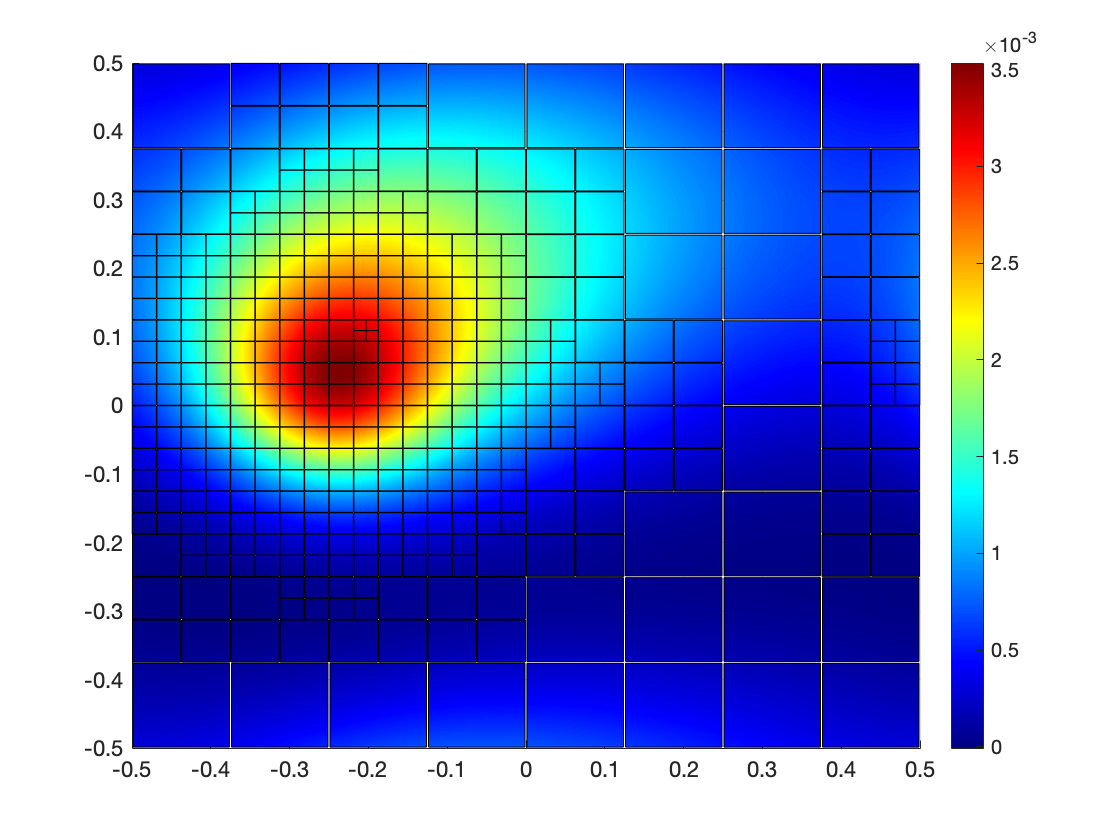}
}
\\
\subfigure[t=0.07, \#leaf boxes=580]{
\includegraphics[width=0.3\textwidth]{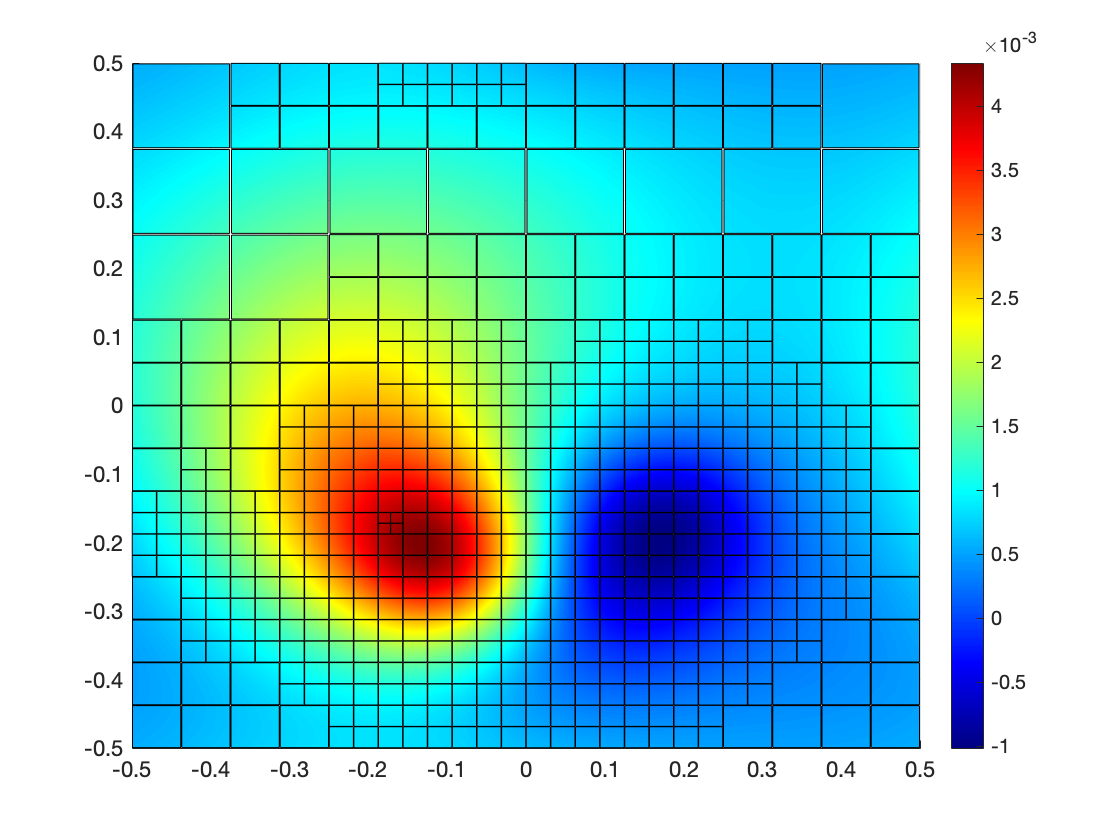}
}
\subfigure[t=0.09, \#leaf boxes=664]{
\includegraphics[width=0.3\textwidth]{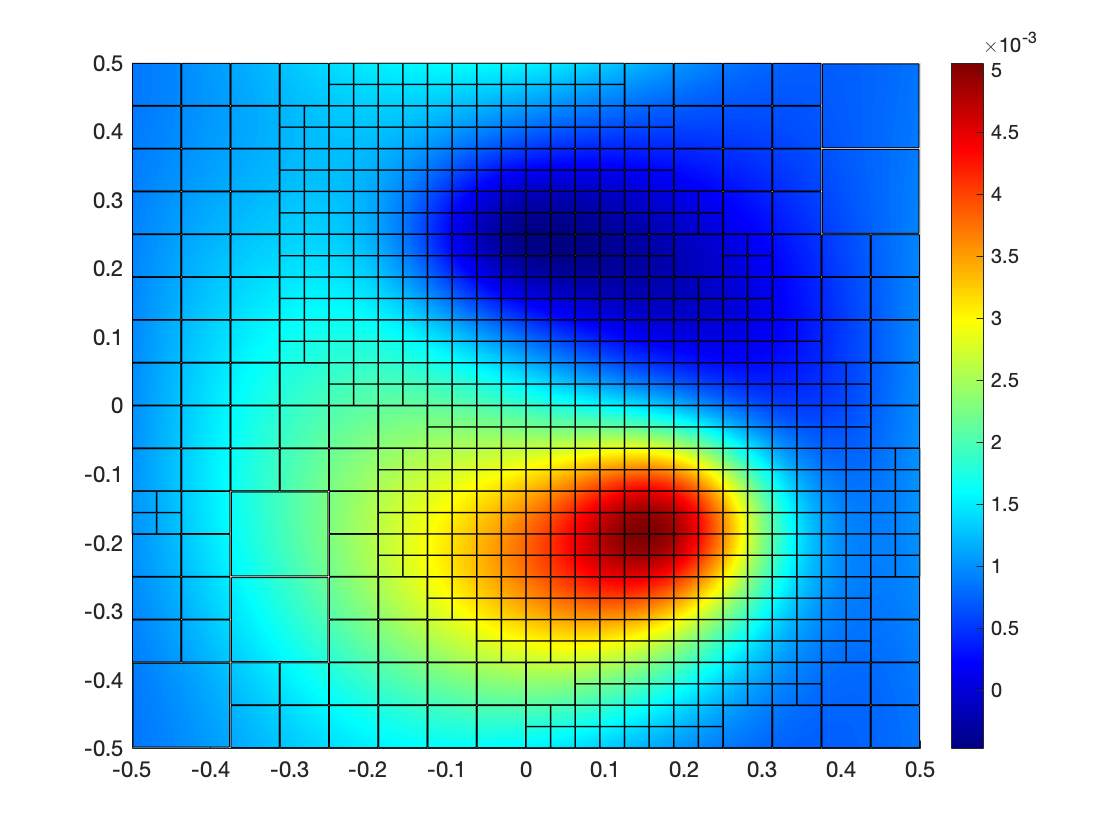}
}
\subfigure[t=0.1, \#leaf boxes=652]{
\includegraphics[width=0.3\textwidth]{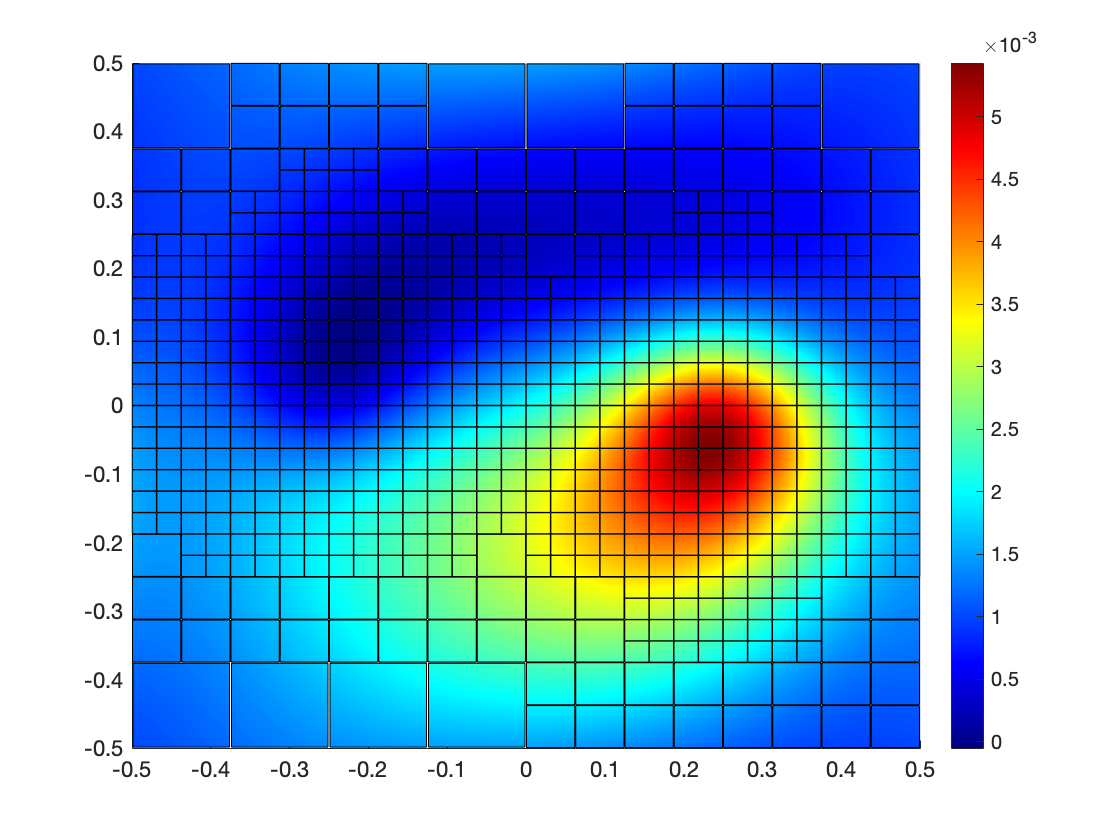}
}
\caption{Evolution of the solution to the inhomogeneous, periodic heat equation with 
forcing function defined in \eqref{heat_forcing}. Each subfigure displays both the 
numerical solution and the corresponding level-restricted quadtree at selected time steps.
 The adaptive trees are generated using our automatic refinement and coarsening strategy, 
and functions on each leaf node are discretized using a scaled $8 \times 8$ 
Chebyshev grid.}
\label{heat_pictures}
\end{figure}

Simulations are performed up to a final time \( T = 0.1 \) with an FGT tolerance 
\( \epsilon = 10^{-9} \). The adaptive spatial grid, which evolves over time to accurately 
capture the features of the solution, is shown in Fig. \ref{heat_pictures} at 
various intermediate times.
Functions on each leaf node in the tree are discretized using a scaled $8 \times 8$ 
Chebyshev grid.
The adaptive tree is generated using the automatic refinement and coarsening strategy
outlined above. Table
\ref{table_heat1} presents the \( L^2 \) norm of the error for 
various time steps using second and fourth order marching schemes. 
The observed convergence rates are consistent with  the theoretical order of accuracy.
Since the scheme is linear scaling and the FGT performance is insensitive to the degree
of adaptivity, a useful rule of thumb is that a time step can be computed in
approximately $(N/10^6)$ seconds on a single core. That is, the combination of
one FGT and one iteration of adaptive refinement and coarsening can be computed at about 
one million points per second. 

\begin{table}[htbp]
 \caption{Convergence data for the forced heat equation (Example 1).
$N_{\rm step}$ denotes the number of time steps and the $L^2$ error
is estimated using the finest grid solution as the reference solution.
$k$ is the estimated order of convergence from the data.
The subscripts $AB2, AB4$ indicate the use of a second or fourth order Adams-Moulton method,
respectively. 
}
  \label{table_heat1}
  \centering
  \begin{tabular}{c|c|c}
    \toprule
$N_{\rm step}$ & $[L^2\ {\rm error}]_{AB2}$ & $k_{AB2}$ \\
\midrule
$2^{10}$& $1.5\times 10^{-5}$   & 1.99  \\
$2^{11}$& $3.8\times 10^{-6}$   &1.99   \\
$2^{12}$& $9.5\times 10^{-7}$   & 1.99  \\
$2^{13}$& $2.4\times 10^{-7}$  &2.00   \\
$2^{14}$& $5.9\times 10^{-8}$   &  \\
    \bottomrule
  \end{tabular}
\qquad \qquad
  \begin{tabular}{c|c|c}
    \toprule
$N_{\rm step}$ & $[L^2\ {\rm error}]_{AB4}$ & $k_{AB4}$ \\
\midrule
$2^6$& $3.7\times 10^{-5}$    & 3.89  \\
$2^7$&  $2.5\times 10^{-6}$   &3.97   \\
$2^8$& $1.6\times 10^{-7}$    & 3.98  \\
$2^9$& $1.0\times 10^{-8}$  &3.80   \\
$2^{10}$&  $7.2\times 10^{-10}$  &   \\

    \bottomrule
  \end{tabular}
\end{table}

\subsection{A reaction-diffusion system}

\begin{figure}[htbp]
\centering
\subfigure[t=10s, \#leaf boxes=946]{
\includegraphics[width=0.3\textwidth]{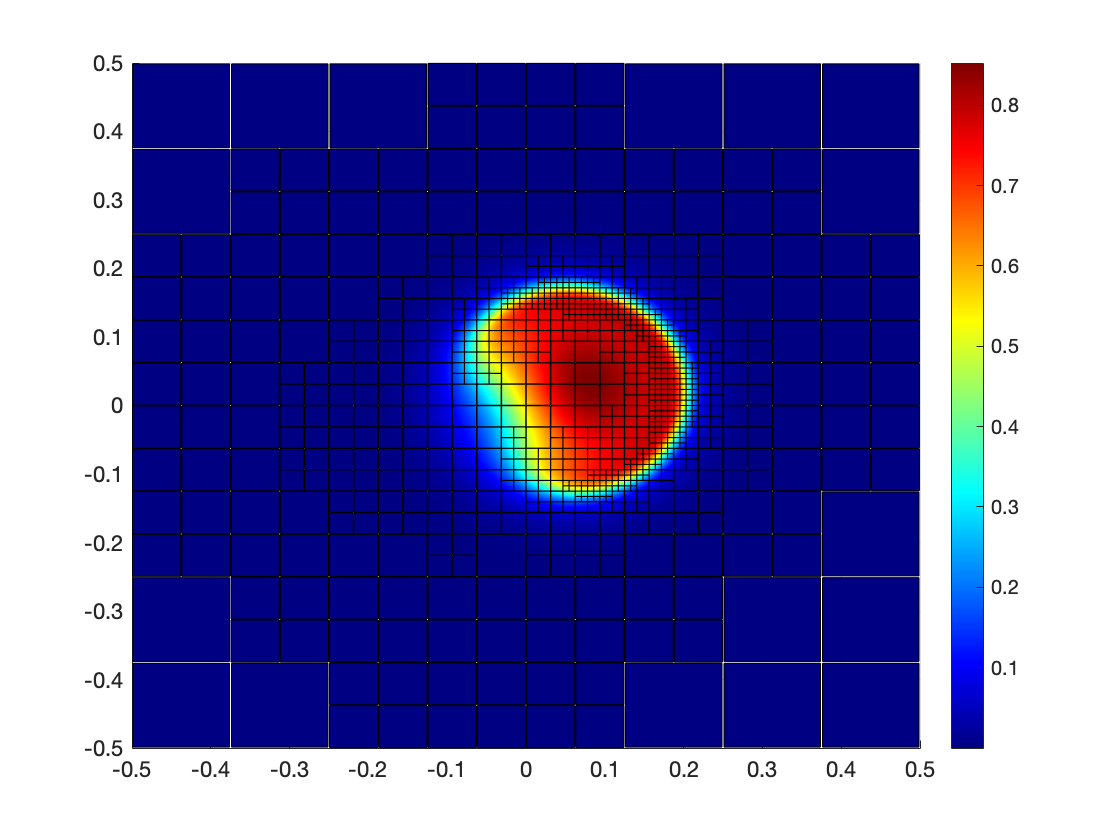}
}
\subfigure[t=100s, \#leaf boxes=1108]{
\includegraphics[width=0.3\textwidth]{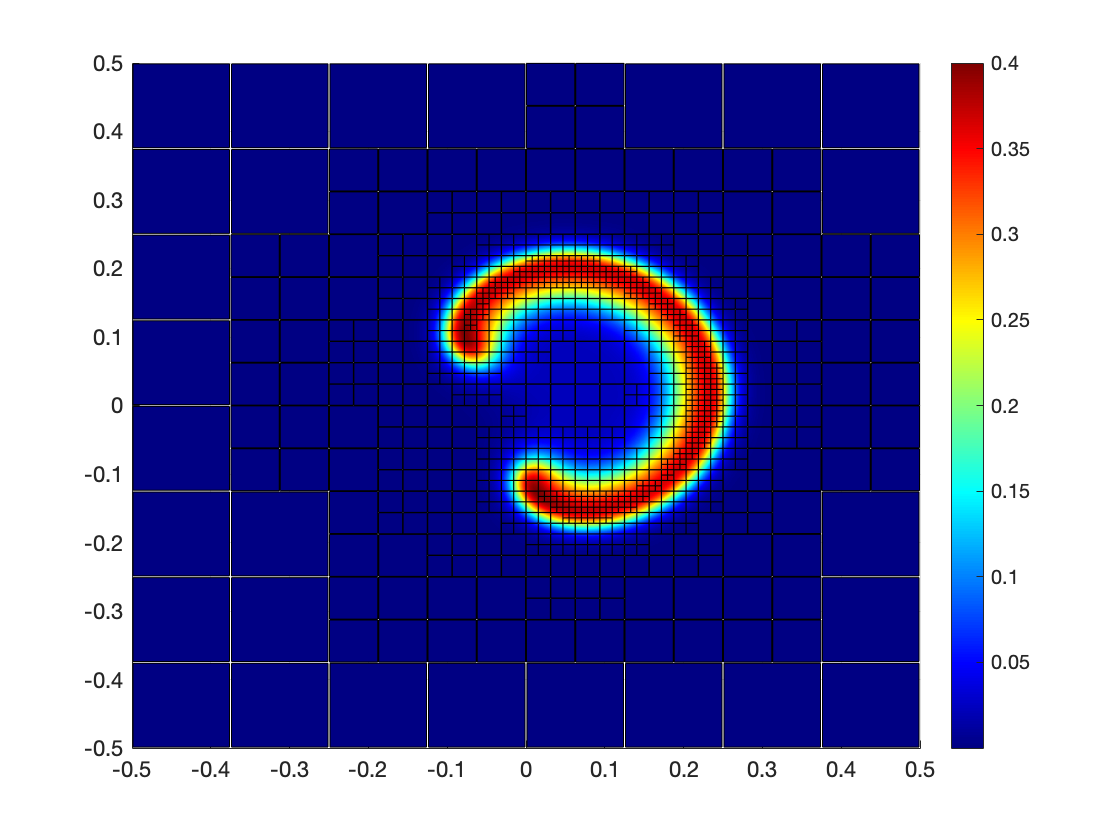}
}
\subfigure[t=400s, \#leaf boxes=1420]{
\includegraphics[width=0.3\textwidth]{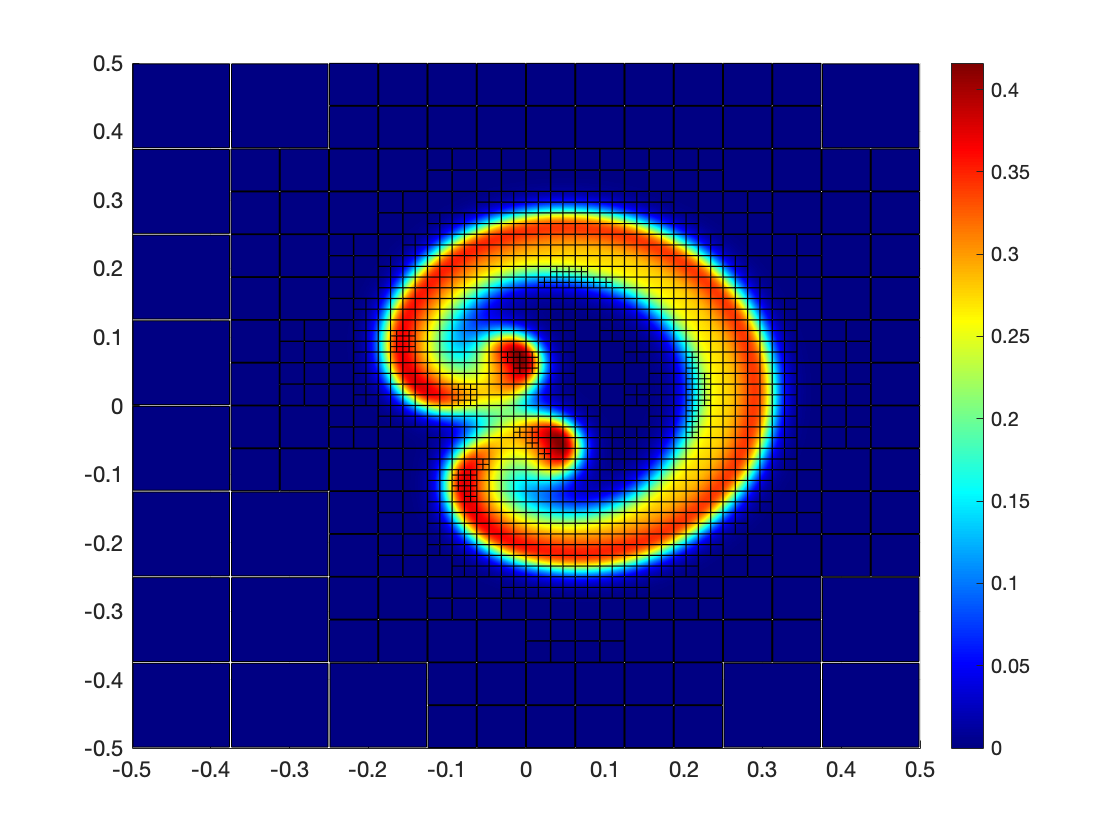}
}
\\
\subfigure[t=600s, \#leaf boxes=1747]{
\includegraphics[width=0.3\textwidth]{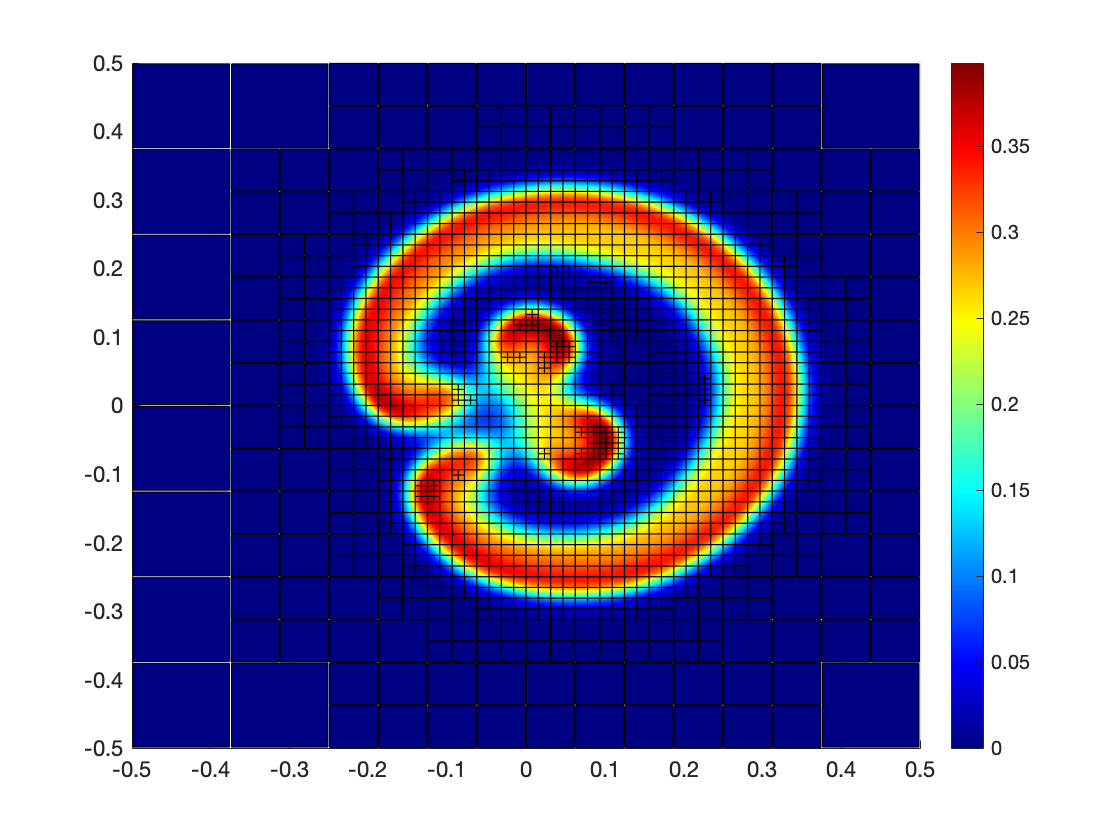}
}
\subfigure[t=800s, \#leaf boxes=2068]{
\includegraphics[width=0.3\textwidth]{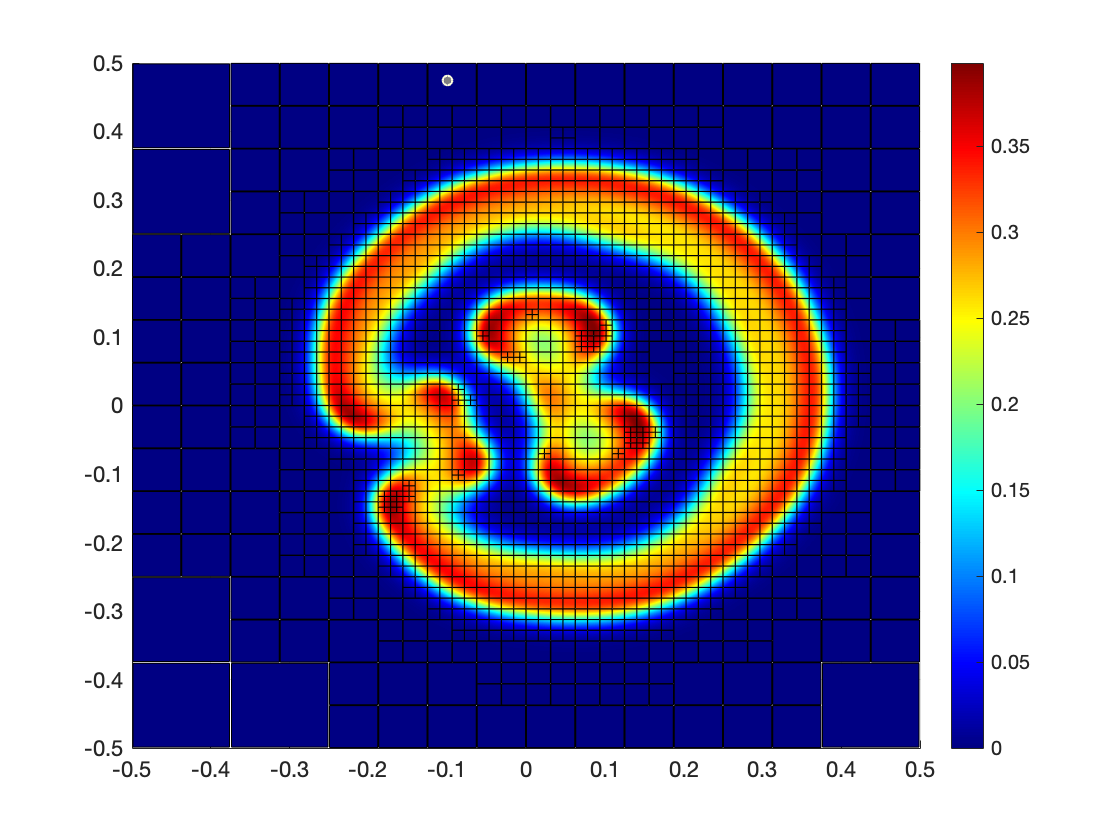}
}
\subfigure[t=1000s, \#leaf boxes=2443]{
\includegraphics[width=0.3\textwidth]{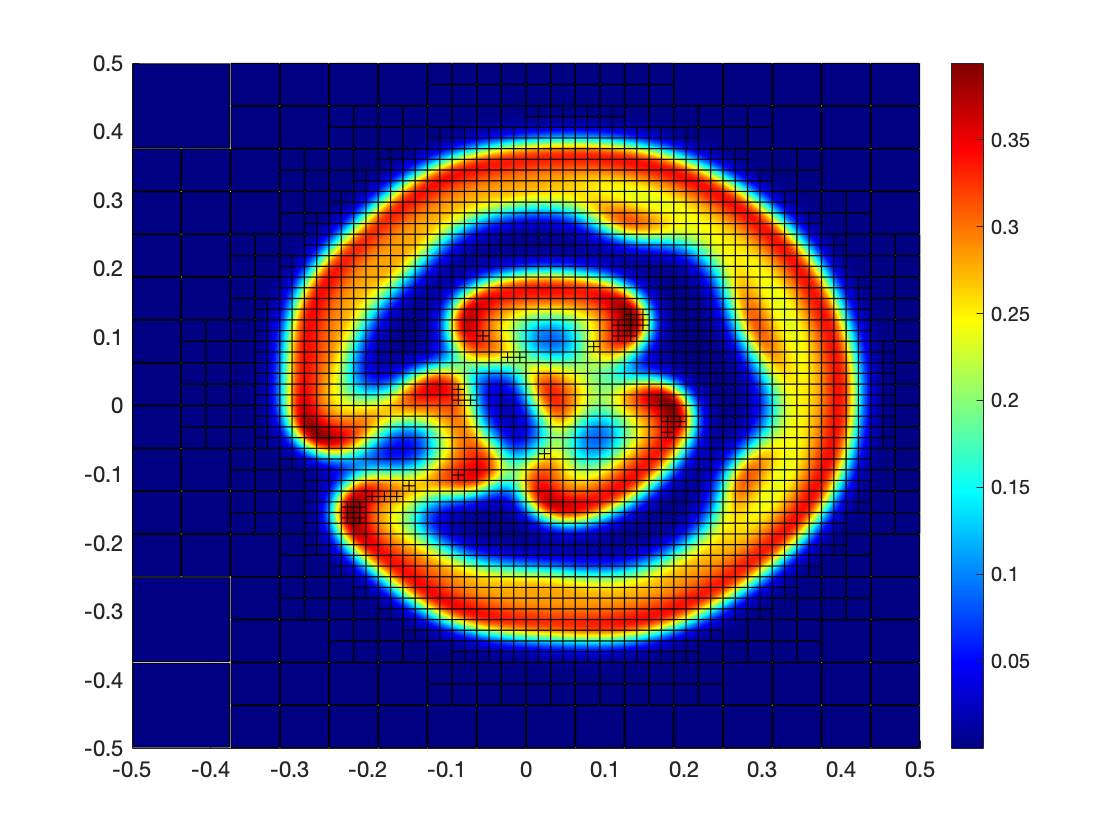}
}
\caption{Evolution of the $v$-component of the Gray-Scott reaction-diffusion 
system at selected time steps. Each panel shows the numerical solution overlaid 
with the corresponding adaptive quadtree mesh. The adaptive meshes are generated 
using our automatic refinement and coarsening strategy
with an error tolerance $\epsilon = 10^{-9}$. The spatial discretization
employs an eighth-order tensor product Chebyshev approximation on each leaf node.
}
\label{heat_sys_pictures}
\end{figure}

For our second example, we consider a reaction-diffusion system using the implicit 
Adams-Moulton integrator discussed in section \ref{sec:AMmethod}. 
In particular, we examine the performance of our scheme on
the Gray-Scott model\cite{gray1990chemical}, a widely studied system that
is known to generate intricate spatio-temporal patterns:
We assume that
\begin{eqnarray}
u_t &= D_u \Delta u - u v^2 + \gamma (1 - u), \nonumber\\
v_t &= D_v \Delta v + u v^2 - (\gamma + \kappa) v,
\end{eqnarray}
with initial conditions 
\begin{equation}
\begin{aligned}
u(x, y, 0) &= 1.0 - \exp\left(-80.0 \left((x + 0.05)^2 + (y + 0.02)^2\right)\right),\\
v(x, y, 0) &= \exp\left(-80.0 \left((x - 0.05)^2 + (y - 0.02)^2\right)\right). 
\end{aligned}
\end{equation}
We use the parameter values \( D_u = 2 \times 10^{-5} \), \( D_v = 1 \times 10^{-5} \), \( \gamma = 0.04 \), and \( \kappa = 0.1 \).
 Simulations are carried out to a final time $T = 1000$ with 
a fixed FGT tolerance $\epsilon = 10^{-9}$. A second-order Adams-Moulton scheme is employed to compute the numerical solution, and the evolution of the variable $v$ is shown
at a few intermediate time steps in Fig. \ref{heat_sys_pictures}, along with the associated time-varying adaptive spatial grid used to resolve the solution.
Convergence results are reported in Table \ref{sys_tab1} for times $T=10$ and $T=200$
using second and fourth order accurate schemes, respectively.

\begin{table}[htbp]
 \caption{Convergence data for the Gray-Scott reaction diffusion system.
$N_{\rm step}$ denotes the number of time steps and the $L^2$ error
is estimated using the finest grid solution as the reference solution.
$k$ is the estimated order of convergence from the data.
The subscripts $AB2, AB4$ indicate the use of a second or fourth order Adams-Moulton method,
respectively. 
}
\label{sys_tab1}
  \centering
  \begin{tabular}{l|l|l}
    \toprule
$N_{\rm step}$ & $[L^2\ {\rm error}]_{AB2}$ & $k_{AB2}$ \\
\midrule
$100$& $1.2 \times 10^{-5}$ &2.0\\
$200$& $3.1 \times 10^{-6}$  &2.0\\
$400$& $7.6 \times 10^{-7}$ &2.0\\
$800$& $1.9 \times 10^{-7}$  & \\
    \bottomrule
  \end{tabular}
\qquad \qquad
\begin{tabular}{l|l|l}
    \toprule
$N_{\rm step}$ & $[L^2\ {\rm error}]_{AB4}$ & $k_{AB4}$ \\
\midrule
& & \\
$300$ & $5.1 \times 10^{-5}$ & 3.7\\
$600$ & $4.0 \times 10^{-6}$  & 4.1\\
$1200$& $2.3 \times 10^{-7}$ & \\
    \bottomrule
  \end{tabular}
\end{table}

A breakdown of the computational cost of
our adaptive solver for the Gray-Scott model is provided in Table \ref{sys_tab4}.
The metric used is {\em grid points processed per second} for each time step
in our single core implementation,
which we denote by ``Rate". 
We show the time and throughput rate for the full time step and for each of the
component modules: the
FGT, the solution of the nonlinear systems at each grid point
using Newton iteration, and adaptive refinement/coarsening. 
The results demonstrate sustained high performance,
over varying problem sizes and levels of adaptivity.

\begin{table}[htbp]
  \caption{Throughput (in grid points per second) for one time step
for different components of the 
adaptive solver applied to the Gray-Scott model.
The subscripts [total, FGT, Newton, adapt] refer to the overall cost, the FGT cost, the 
cost of solving the nonlinear systems at each grid point, and the cost of 
refinement/coarsening, respectively. $N_{\text{pts}} = N_{\rm leaf} \times 64$
is the total number of points in the adaptive grid.
}
  \centering
  \begin{tabular}{r|r|l|l|l|l}
    \toprule
    $N_{\text{step}}$ & $N_{\text{pts}}$ & $T_{\text{total}}$ & Rate & $T_{\text{FGT}}$ & Rate \\
    \midrule
    10     & 35904   & $9.2 \times 10^{-2}$  & $3.9 \times 10^5$ & $9.4 \times 10^{-3}$ & $3.8 \times 10^6$ \\
    500    & 97856  & $2.8 \times 10^{-1}$ & $3.5 \times 10^5$ & $7.4 \times 10^{-2}$ & $1.3 \times 10^6$ \\
    1000   & 94528  & $2.6 \times 10^{-1}$ & $3.7 \times 10^5$ & $5.1 \times 10^{-2}$ & $1.9 \times 10^6$ \\
    3000   & 111680  & $2.8 \times 10^{-1}$ & $4.0 \times 10^5$ & $4.9 \times 10^{-2}$ & $2.3 \times 10^6$ \\
    5000   & 139328  & $3.2 \times 10^{-1}$ & $4.3 \times 10^5$ & $5.5 \times 10^{-2}$ & $2.5 \times 10^6$ \\
    8000   & 176448  & $3.9 \times 10^{-1}$ & $4.5 \times 10^5$ & $6.7 \times 10^{-2}$ & $2.7 \times 10^6$ \\
    10000  & 208448  & $4.6 \times 10^{-1}$ & $4.5 \times 10^5$ & $7.4 \times 10^{-2}$ & $2.8 \times 10^6$ \\
    \midrule
    $N_{\text{step}}$ & $N_{\text{leaf}}$ & $T_{\text{Newton}}$ & Rate & $T_{\text{adapt}}$ & Rate \\
    \midrule
    10     & 35904   & $4.0 \times 10^{-2}$ & $9.0 \times 10^5$  & $2.9 \times 10^{-2}$ & $1.2 \times 10^6$ \\
    500    & 97856  & $8.7 \times 10^{-2}$ & $1.1 \times 10^6$  & $6.1 \times 10^{-2}$ & $1.6 \times 10^6$ \\
    1000   & 94528  & $8.1 \times 10^{-2}$ & $1.2 \times 10^6$  & $5.7 \times 10^{-2}$ & $1.7 \times 10^6$ \\
    3000   & 111680  & $1.0 \times 10^{-1}$ & $1.1 \times 10^6$  & $5.7 \times 10^{-2}$ & $2.0 \times 10^6$ \\
    5000   & 139328  & $1.2 \times 10^{-1}$ & $1.1 \times 10^6$  & $6.7 \times 10^{-2}$ & $2.1 \times 10^6$ \\
    8000   & 176448  & $1.5 \times 10^{-1}$ & $1.2 \times 10^6$  & $8.8 \times 10^{-2}$ & $2.0 \times 10^6$ \\
    10000  & 208448  & $1.8 \times 10^{-1}$ & $1.2 \times 10^6$  & $9.8 \times 10^{-2}$ & $2.1 \times 10^6$ \\
    \bottomrule
  \end{tabular}

  \label{sys_tab4}
\end{table}

\subsection{Unsteady Stokes flow}

In our third example, we solve the unsteady Stokes equations with exact solution
\begin{equation}
\begin{aligned}
u(x, y,t) &= \pi e^{\sin(2\pi x)} e^{\sin(2\pi y)} \cos(2\pi y) \sin^2(t),\\
v(x, y,t) &= -\pi e^{\sin(2\pi x)} e^{\sin(2\pi y)} \cos(2\pi x) \sin^2(t), \\
p(x, y,t) &= e^{\cos(2\pi x) \sin(2\pi y)} \sin^2(t).
\end{aligned}
\end{equation}
We use an error tolerance $\epsilon = 10^{-9}$ for function resolution, with convergence
data presented in Table \ref{tab:convergence-unsteady}.

\begin{table}[htbp]
 \caption{Convergence data for the unsteady Stokes equations.
$N_{\rm step}$ denotes the number of time steps and the $L^2$ error
is estimated using the finest grid solution as the reference solution.
$k$ is the estimated order of convergence from the data.
The subscripts $AB2, AB4$ indicate the use of a second or fourth order Adams-Moulton method,
respectively. 
The second order code is run to a final time $T=0.1$ (left) and 
the fourth order code is run to a final time $T=0.1$ (right).}
  \centering
  \begin{tabular}{c|c|c}
    \toprule
$N_{\rm step}$ & $[L^2\ {\rm error}]_{AB2}$ & $k_{AB2}$ \\
\midrule
$2^{7}$& $1.6\times 10^{-3}$   & 2.0  \\
$2^{8}$& $3.9\times 10^{-4}$   & 2.0  \\
$2^{9}$& $9.6\times 10^{-5}$   & 2.0  \\
$2^{10}$& $2.5\times 10^{-5}$   & 2.0  \\
$2^{11}$& $6.2\times 10^{-6}$   &  \\
    \bottomrule
  \end{tabular}
\qquad \qquad
  \begin{tabular}{c|c|c}
    \toprule
$N_{\rm step}$ & $[L^2\ {\rm error}]_{AB4}$ & $k_{AB4}$ \\
\midrule
$2^5$& $1.5\times 10^{-1}$    & 2.5  \\
$2^6$& $2.7\times 10^{-2}$    & 3.0  \\
$2^7$&  $3.5\times 10^{-3}$   & 3.5   \\
$2^8$& $3.1\times 10^{-4}$    & 3.8  \\
$2^9$& $2.2\times 10^{-5}$    & 3.7   \\
$2^{10}$&  $1.7\times 10^{-6}$  & 3.9  \\
$2^{11}$&  $1.1\times 10^{-7}$  & 3.9  \\
$2^{12}$&  $7.6\times 10^{-9}$  &   \\
    \bottomrule
  \end{tabular}
  \label{tab:convergence-unsteady}
\end{table}

\begin{figure}[htbp]
\centering
\subfigure[t=0.01s \#leaf boxes=1882]{
\includegraphics[width=0.3\textwidth]{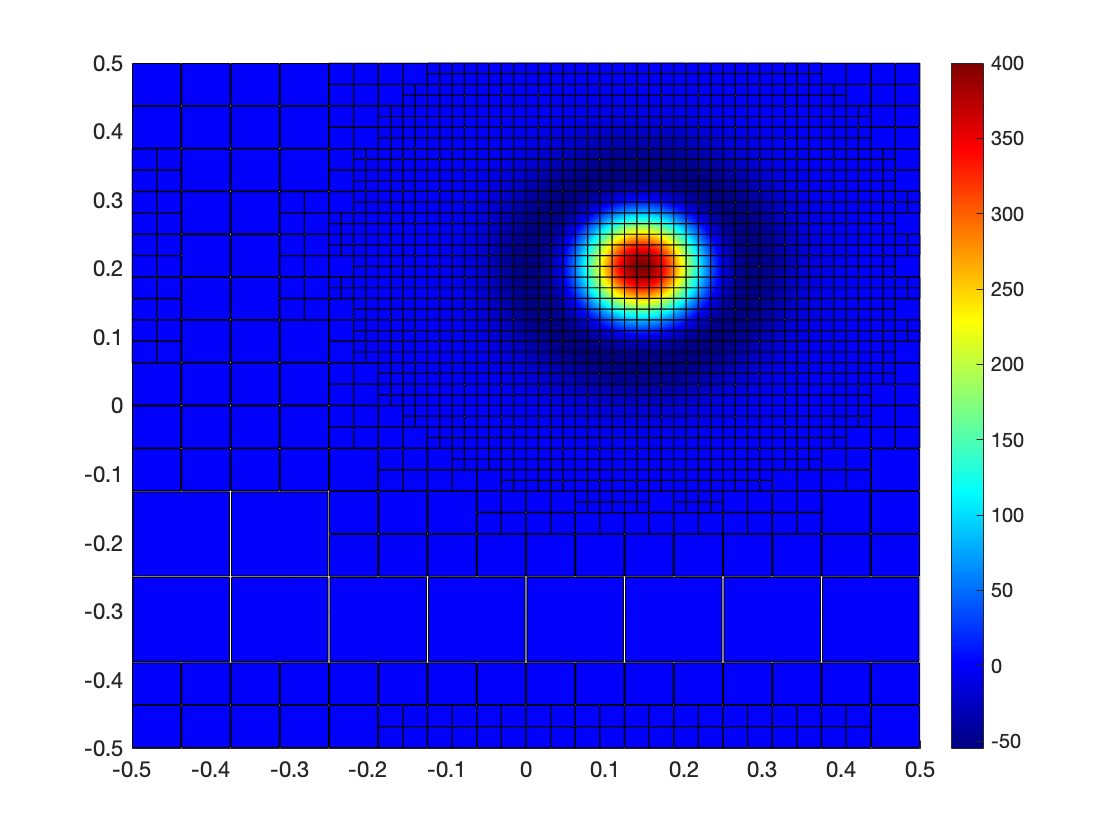}
}
\subfigure[t=0.03s \#leaf boxes=1816]{
\includegraphics[width=0.3\textwidth]{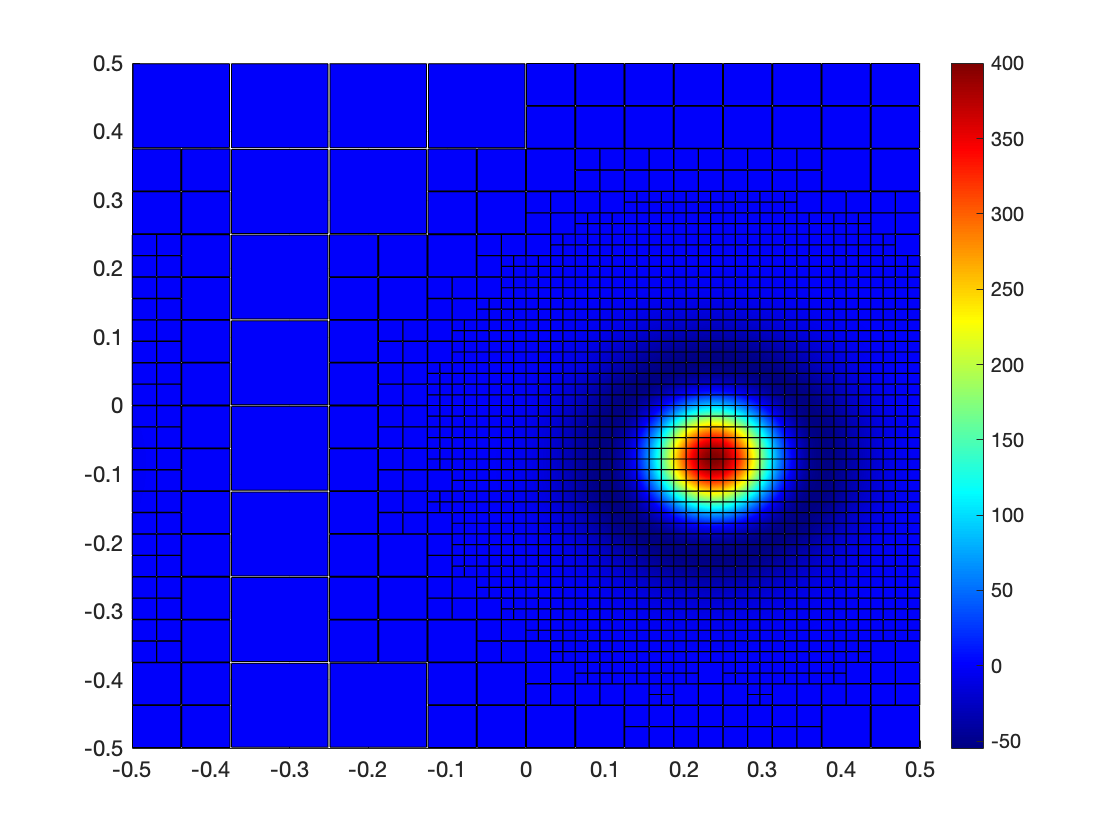}
}
\subfigure[t=0.05s \#leaf boxes=1854]{
\includegraphics[width=0.3\textwidth]{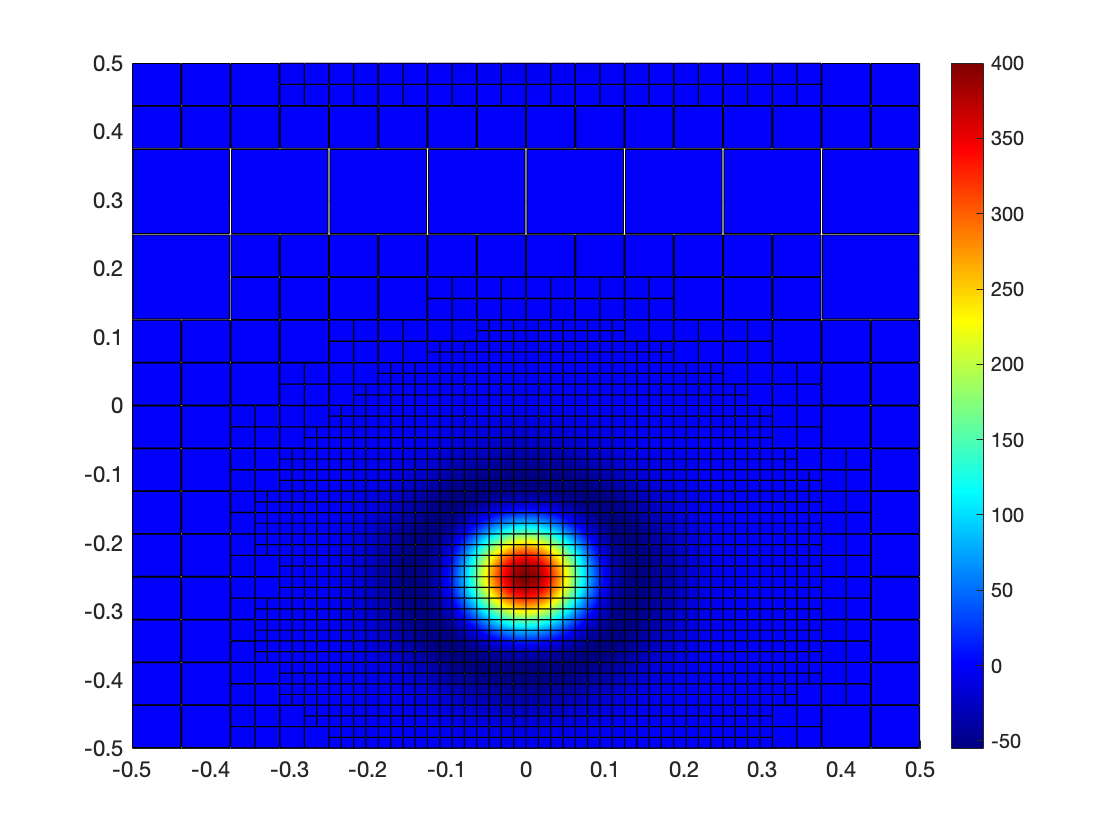}
}
\\
\subfigure[t=0.07s \#leaf boxes=1852]{
\includegraphics[width=0.3\textwidth]{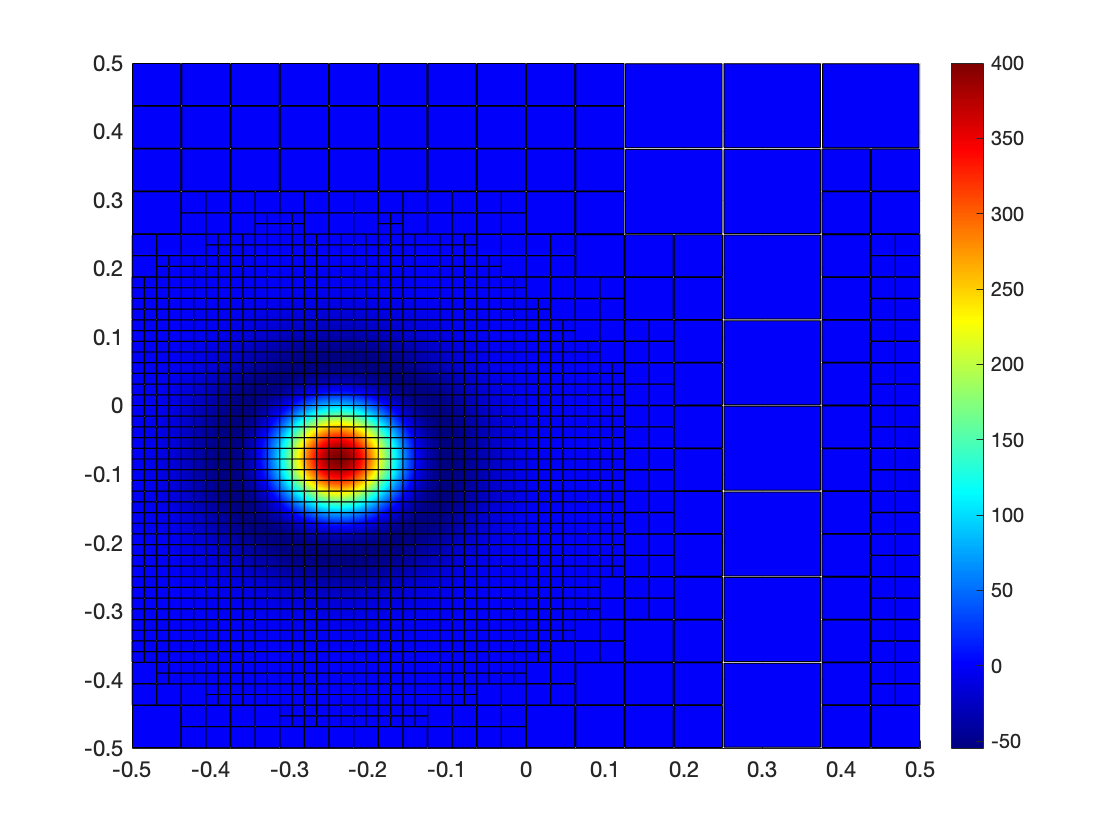}
}
\subfigure[t=0.09s \#leaf boxes=1912]{
\includegraphics[width=0.3\textwidth]{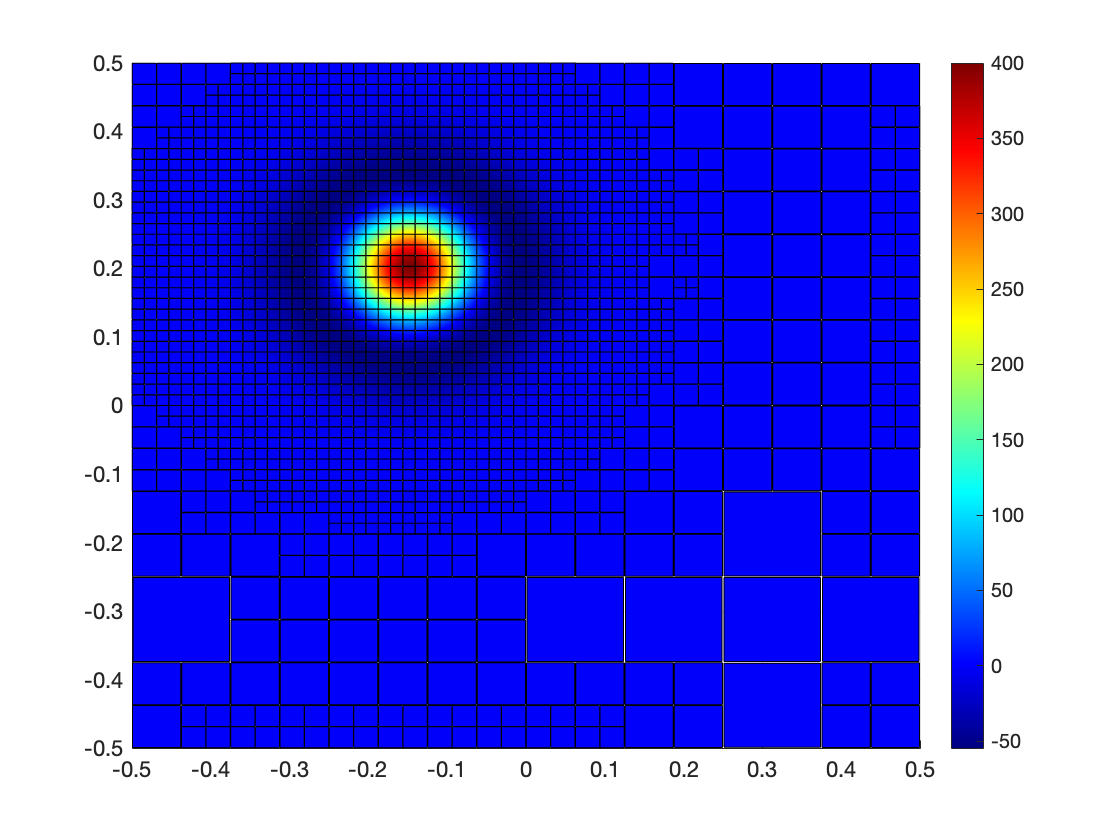}
}
\subfigure[t=0.1s \#leaf boxes=1804]{
\includegraphics[width=0.3\textwidth]{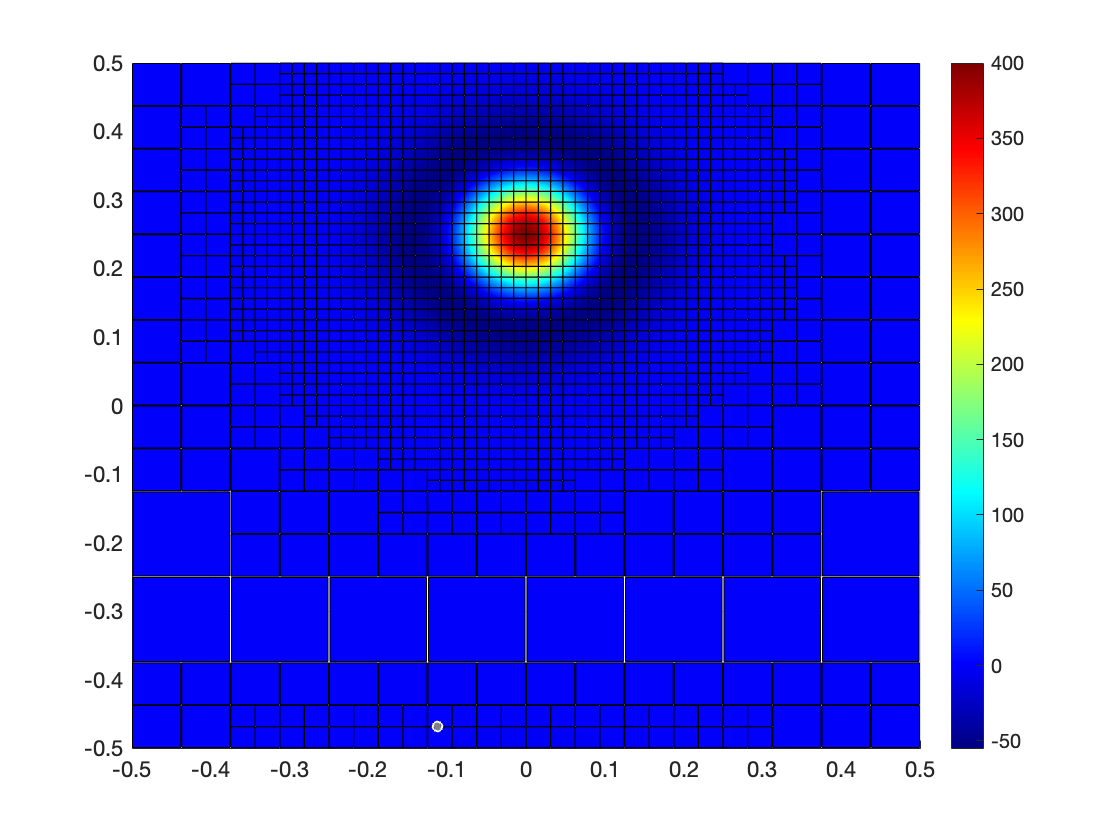}
}
\caption{Evolution of the vorticity for the  unsteady Stokes equations\eqref{eq-us2} 
on $D = [-0.5, 0.5]^2$ at selected time steps up to the final time $T=0.1$. 
Each subfigure shows the numerical solution overlaid with the corresponding adaptive quadtree mesh.} 
\label{us_pictures}
\end{figure}

In order to profile the performance of our adaptive solver, we consider
the unsteady Stokes equations with an explicit fourth order Adams Bashforth method 
and an exact solution given by
\begin{equation}
\begin{aligned}
u(x,y,t) &= \sum_{j\in \mathbf{Z}} \frac{2(y - y_j)}{\delta} e^{-\frac{(x - x_j)^2 + (y - y_j)^2}{\delta}},\\
v(x,y,t) &= \sum_{j\in \mathbf{Z}} \frac{2(x - x_j)}{\delta} e^{-\frac{(x - x_j)^2 + (y - y_j)^2}{\delta}},
\end{aligned}\label{eq-us2}
\end{equation}
with $x_j = 0.25 \sin(20\pi t) + j$
and $y_j = 0.25 \cos(20\pi t) + j$.
In Fig. \ref{us_pictures}, we show the vorticity contours at 
selected time steps up to the final time $T = 0.1$.

As above, we use the metric of grid points processed per second for each time step on 
a single core, with $\dt=0.01$. 
In Table \ref{us_tab}, we show the time and throughput rate for the full time step 
and for each of the component modules: the FGT, the Helmholtz 
projection using the FMM, and
adaptive refinement/coarsening. 
Note that the throughput for a single time step is about $180,000$ points per second.

\begin{table}[htbp]
\caption{Throughput (in grid points per second) for one time step
of the adaptive solver applied to the unsteady Stokes equations.
The subscripts [total, FGT, FMM, adapt] refer to the overall cost, the FGT cost, the 
cost of Helmholtz projection using the FMM, and the cost of 
refinement/coarsening, respectively. $N_{\text{pts}} = N_{\rm leaf} \times 64$
is the total number of points in the adaptive grid.
}
\centering
\begin{tabular}{l|l|l|l|l|l}
\toprule
$N_{\text{total}}$ & $N_{\text{leaf}}$ & $T_{\text{step}}$ & Rate & $T_{\text{FGT}}$ & Rate \\
\midrule
6 & 115,648 & $6.3 \times 10^{-1}$ & $1.8 \times 10^5$ & $5.6 \times 10^{-2}$ & $2.1 \times 10^6$ \\
15 & 121,984 & $6.6 \times 10^{-1}$ & $1.9 \times 10^5$ & $4.0 \times 10^{-2}$ & $3.0 \times 10^6$ \\
35 & 120,832 & $6.6 \times 10^{-1}$ & $1.8 \times 10^5$ & $4.0 \times 10^{-2}$ & $3.0 \times 10^6$ \\
50 & 115,456 & $6.4 \times 10^{-1}$ & $1.8 \times 10^5$ & $3.9 \times 10^{-2}$ & $2.9 \times 10^6$ \\
78 & 116,416 & $6.4 \times 10^{-1}$ & $1.8 \times 10^5$ & $4.0 \times 10^{-2}$ & $2.9 \times 10^6$ \\
96 & 116,416 & $6.4 \times 10^{-1}$ & $1.8 \times 10^5$ & $3.8 \times 10^{-2}$ & $3.1 \times 10^6$ \\
\midrule
$N_{\text{total}}$ & $N_{\text{leaf}}$ & $T_{\text{FMM}}$ & Rate & $T_{\text{adapt}}$ & Rate \\
\midrule
6 & 115,648 & $7.1 \times 10^{-2}$ & $2.2 \times 10^6$ & $2.6 \times 10^{-2}$ & $4.5 \times 10^6$ \\
15 & 121,984 & $7.8 \times 10^{-2}$ & $2.1 \times 10^6$ & $3.1 \times 10^{-2}$ & $4.0 \times 10^6$ \\
35 & 120,832 & $7.7 \times 10^{-2}$ & $2.1 \times 10^6$ & $3.0 \times 10^{-2}$ & $4.0 \times 10^6$ \\
50 & 115,456 & $7.9 \times 10^{-2}$ & $2.0 \times 10^6$ & $3.1 \times 10^{-2}$ & $3.7 \times 10^6$ \\
78 & 116,416 & $7.0 \times 10^{-2}$ & $2.2 \times 10^6$ & $2.8 \times 10^{-2}$ & $4.2 \times 10^6$ \\
96 & 116,416 & $7.5 \times 10^{-2}$ & $2.1 \times 10^6$ & $3.0 \times 10^{-2}$ & $3.8 \times 10^6$ \\
\bottomrule
\end{tabular}
\label{us_tab}
\end{table}

\subsection{The Navier-Stokes equations}

\begin{figure}[htbp]
\centering
\subfigure[t=0.008s \#leaf boxes=1600]{
\includegraphics[width=0.3\textwidth]{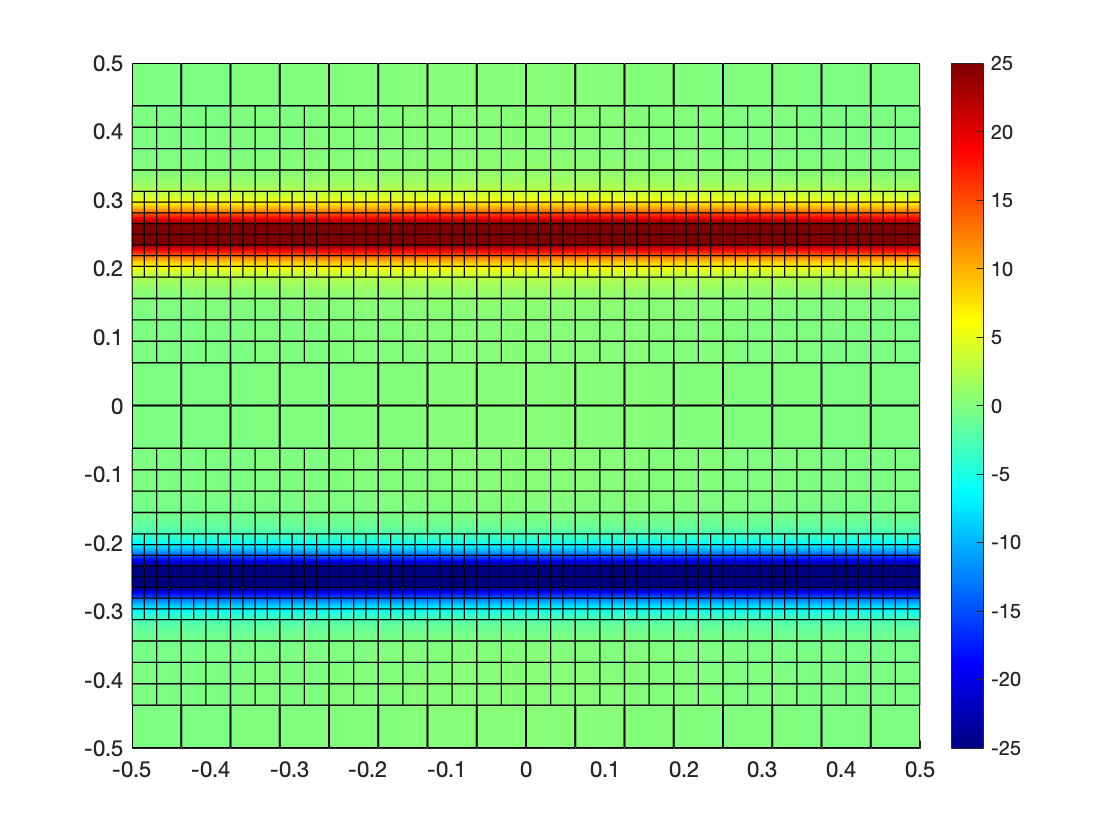}
}
\subfigure[t=0.24s \#leaf boxes=1792]{
\includegraphics[width=0.3\textwidth]{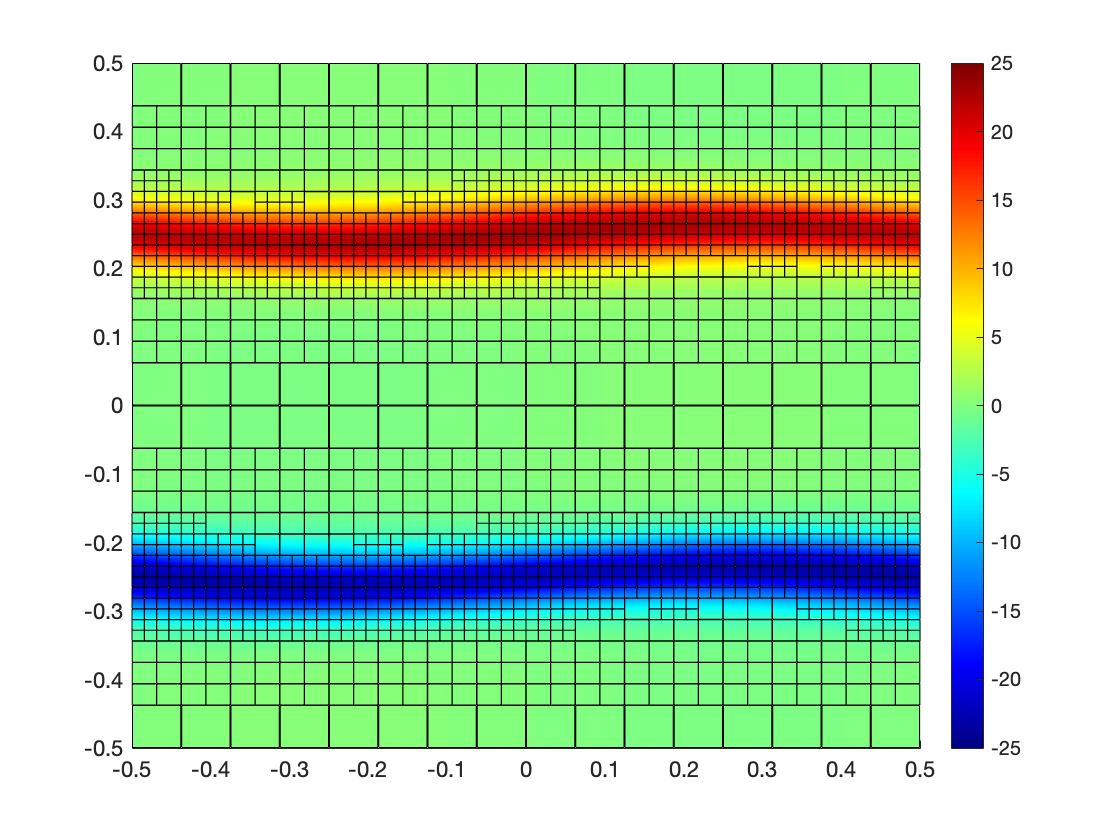}
}
\subfigure[t=0.6s \#leaf boxes=1696]{
\includegraphics[width=0.3\textwidth]{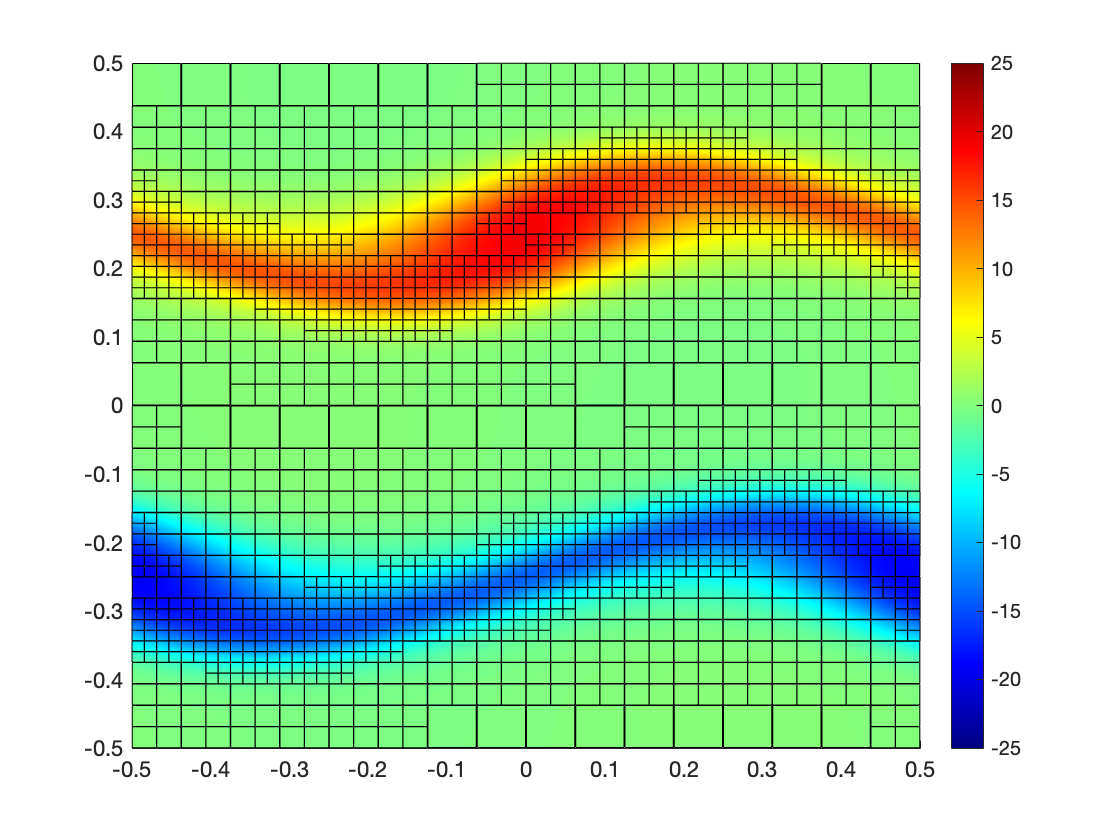}
}
\\
\subfigure[t=0.96s \#leaf boxes=1864]{
\includegraphics[width=0.3\textwidth]{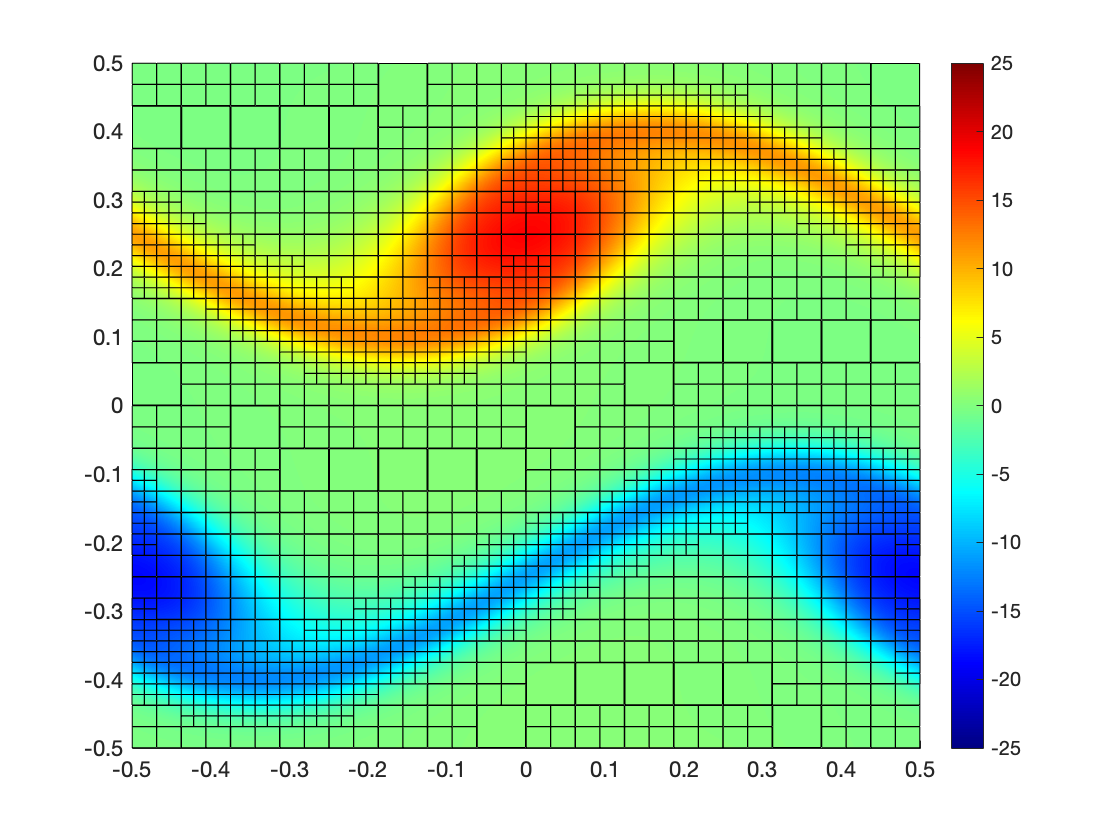}
}
\subfigure[t=1.08s \#leaf boxes=1984]{
\includegraphics[width=0.3\textwidth]{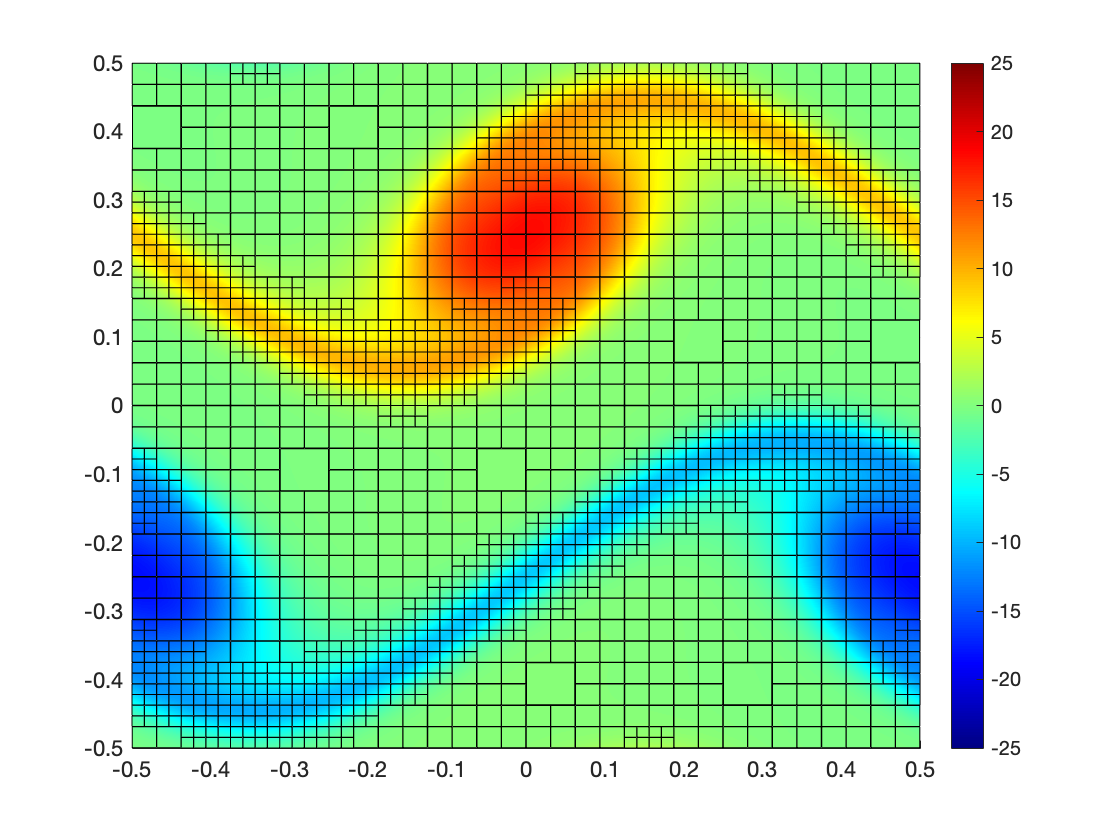}
}
\subfigure[t=1.2s \#leaf boxes=2152]{
\includegraphics[width=0.3\textwidth]{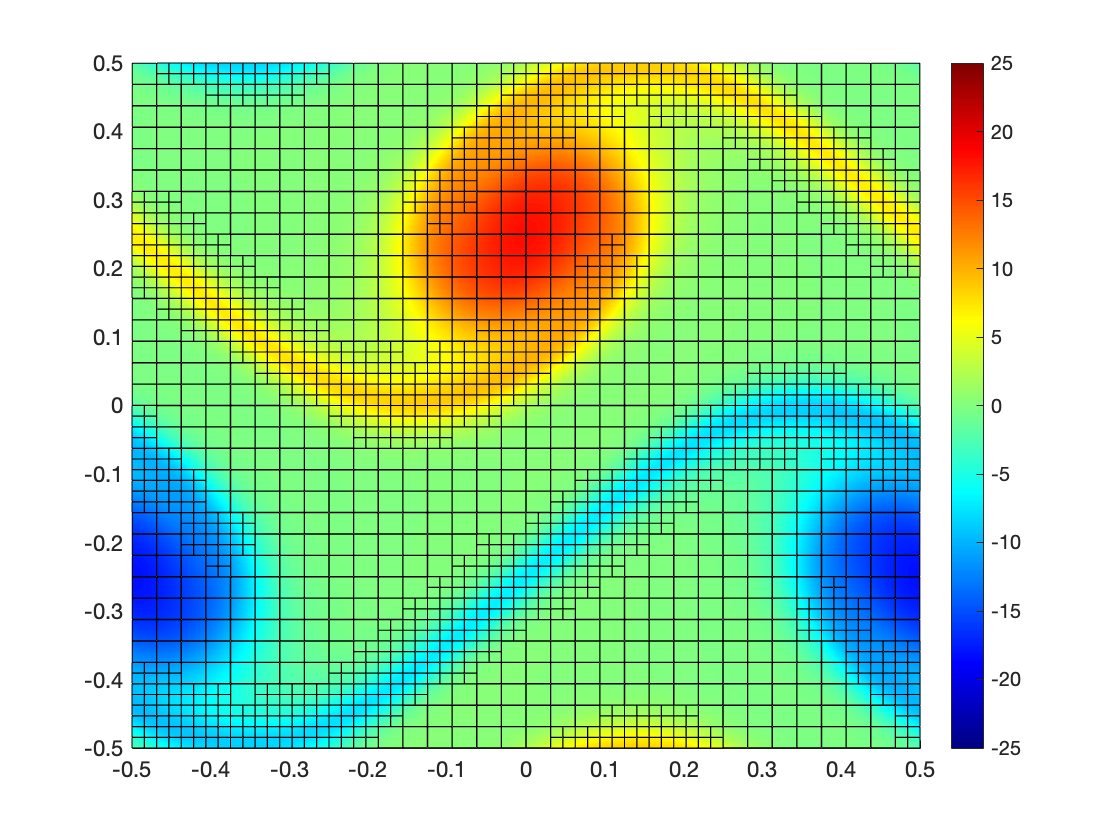}
}
\caption{Evolution of the vorticity field for the double shear layer problem\eqref{double sheer layer_eq} at selected time steps with $\Delta t=0.0008s$. Each subfigure shows the numerical solution overlaid with the corresponding adaptive quadtree.}
\label{ns_pictures}
\end{figure}

In our final example, we consider the classical double shear layer problem \cite{bell1989second,brown1995performance,di2005moving,huang2021stability} governed by the Navier-Stokes equations.
The initial condition is specified as:
\begin{equation}
\begin{aligned}
    u_1(x,y,0)&=
\begin{cases}
\tanh(\rho(y+0.25))& \text{$x \le 0$}\\
\tanh(\rho(-y+0.25))& \text{$x > 0$}
\end{cases} \\
u_2(x,y,0)&= \delta \sin(2\pi x),
\end{aligned}\label{double sheer layer_eq}
\end{equation}
where $\rho$ determines the slope of the shear layer and $\delta$ represents the size of 
the perturbation. In our simulations, we fix $\delta = 0.05$ and $\rho = 30$.
We present snapshots of the vorticity at selected time steps in Fig. \ref{ns_pictures},
with the viscosity set to \( \nu = 0.001 \). As the simulation progresses, the shear layers 
roll up into vortical structures and the adaptive spatial discretization automatically
adjusts the grid to capture the evolving vorticity, maintaining both efficiency and
high accuracy.

As an initial test for our solver, we set 
the viscosity to \( \nu = 0.01 \) and the final time to \( T = 0.4 \).
The $L^2$ error in the solution is assessed using a twice-refined grid with 80 
times steps as a reference solution.
As can be seen in Table \ref{ns_tab1}, the observed convergence rate is consistent with 
the requested fourth-order accuracy of the predictor-corrector scheme.

\begin{table}[htbp]
 \caption{Convergence of a fourth order predictor-corrector scheme
for the double shear layer problem governed by the Navier-Stokes equations.
$k$ is the estimated order of convergence from the data.}

\centering
\begin{tabular}{l|l|l}
\toprule
$N_{\rm step}$ &$L^2$ error & $k$\\
\midrule
10& $1.5 \times 10^{-6}$  &3.3\\
20 & $1.4 \times 10^{-7} $  & 3.7\\
40 & $1.1 \times 10^{-8}$    & 4.0 \\
\bottomrule
\end{tabular}
\label{ns_tab1}
\end{table}

We profile the performance of our solver in Table \ref{sys_tab5}. 
Speed is reported in units of grid points processed per second per time step. 
We show the time and throughput rate for the full time step 
and for each of the component modules: the FGT, the Helmholtz 
projection using the FMM, and adaptive refinement/coarsening. 
The results indicate consistent high performance across each algorithmic component.

\begin{table}[htbp]
\caption{Throughput (in grid points per second) for the adaptive Navier-Stokes solver
applied to the double shear layer problem.
$N_{\rm step}$ denotes the number of time steps.
The subscripts [total, FGT, FMM, adapt] refer to the overall cost, the FGT cost, the 
cost of Helmholtz projection using the FMM, and the cost of 
refinement/coarsening, respectively. $N_{\text{pts}} = N_{\rm leaf} \times 64$
is the total number of points in the adaptive grid at the final time.}
  \centering
  \begin{tabular}{r|r|l|l|l|l}
    \toprule
    $N_{\text{step}}$ & $N_{\text{pts}}$ & $T_{\text{step}}$ & Rate & $T_{\text{FGT}}$ & Rate \\
    \midrule
    $10$   & $102,400$ & $9.0 \times 10^{-1}$ & $1.1 \times 10^{5}$ & $1.6 \times 10^{-2}$ & $2.0 \times 10^{7}$ \\
    $100$  & $110,848$ & $9.6 \times 10^{-1}$ & $1.2 \times 10^{5}$ & $1.8 \times 10^{-2}$ & $1.9 \times 10^{7}$ \\
    $500$  & $116,992$ & $9.5\times 10^{-1}$ & $1.2 \times 10^{5}$ & $1.8 \times 10^{-2}$ & $2.1 \times 10^{7}$ \\
    $700$  & $100,864$ & $9.3\times 10^{-1}$ & $1.1 \times 10^{5}$ & $1.7 \times 10^{-2}$ & $1.8 \times 10^{7}$ \\
    $900$  & $92,416$ & $8.6\times 10^{-1}$ & $1.1 \times 10^{5}$ & $1.5 \times 10^{-2}$ & $1.9 \times 10^{7}$ \\
    $1200$ & $73,216$ & $7.3\times 10^{-1}$ & $1.0 \times 10^{5}$ & $1.2 \times 10^{-2}$ & $1.7 \times 10^{7}$ \\
    $1500$ & $67,072$ & $7.0\times 10^{-1}$ & $9.5 \times 10^{4}$ & $1.1 \times 10^{-2}$ & $1.6 \times 10^{7}$ \\
    \midrule
    $N_{\text{step}}$ & $N_{\text{pts}}$ & $T_{\text{FMM}}$ & Rate & $T_{\text{adapt}}$ & Rate \\
    \midrule
    $10$   & $102,400$ & $7.1 \times 10^{-2}$ & $1.9 \times 10^{6}$ & $8.8 \times 10^{-2}$ & $1.2 \times 10^{6}$ \\
    $100$  & $110,848$ & $7.5 \times 10^{-2}$ & $2.1 \times 10^{6}$ & $9.2 \times 10^{-2}$ & $1.2 \times 10^{6}$ \\
    $500$  & $116,992$ & $7.4 \times 10^{-2}$ & $2.1 \times 10^{6}$ & $8.9 \times 10^{-2}$ & $1.3 \times 10^{6}$ \\
    $700$  & $100,864$ & $7.3 \times 10^{-2}$ & $1.9 \times 10^{6}$ & $8.5 \times 10^{-2}$ & $1.2 \times 10^{6}$ \\
    $900$  & $92,416$ & $6.8 \times 10^{-2}$ & $1.8 \times 10^{6}$ & $8.3 \times 10^{-2}$ & $1.1 \times 10^{6}$ \\
    $1200$ & $73,216$ & $5.6 \times 10^{-2}$ & $1.7 \times 10^{6}$ & $7.8 \times 10^{-2}$ & $9.4 \times 10^{5}$ \\
    $1500$ & $67,072$ & $5.4 \times 10^{-2}$ & $1.7 \times 10^{6}$ & $7.9 \times 10^{-2}$ & $8.5 \times 10^{5}$ \\
    \bottomrule
  \end{tabular}

  \label{sys_tab5}
\end{table}

\clearpage

\section{Conclusions}\label{sec:conclusion}
We have presented a unified approach to parabolic
evolution equations including reaction-diffusion systems,
the unsteady Stokes equations and the full, incompressible
Navier-Stokes equations in a periodic box.
In the linear setting,
these tools permit the creation of stable, explicit, fully automatic,
high-order, space-time adaptive solvers.
In the semi-linear setting, a fully implicit formulation needed 
to handle stiff forcing terms requires only
the solution of uncoupled scalar nonlinear equations at each grid point.
The principal numerical algorithm used here is the (continuous) fast 
Gauss transform. For incompressible flow problems, we also make use of the FMM
to compute the Helmholtz decomposition.
Using heat (and harmonic) layer potentials, one can impose
boundary conditions in complicated interior or exterior domains with comparable
performance in terms of work per grid point and we are currently working on such extensions.

In the present paper, we have assumed that the diffusion coefficient
itself is constant. Developing integral equation methods for piecewise-constant
and piecewise-smooth diffusion coefficients is more involved but also possible,
and a topic of active research. 
Finally, while our current code 
allows for an adaptive grid that evolves with time, the time step itself is assumed to be global. Developing variants that permit local time stepping
is also a topic of active research.

Code libraries for the schemes presented here are being prepared for
release as open source software.

\section*{Acknowledgments}
We would like to thank Travis Askham, Alex Barnett, and Manas Rachh for many
helpful discussions.
Jun Wang and Jie Su were supported in part by the National Natural Science Foundation of China under grant 12301515, and in part by  National Key Research and Development Program of China under grant 2023YFA1008902.

\end{document}